\documentclass[%
 reprint,
 amsmath,amssymb,
 aps,
]{revtex4-2}

\AtBeginDocument{%
  
}

\usepackage{graphicx} 

\makeatletter
\def\l@section{\@tocline{1}{0pt}{0pc}{}{}}
\def\l@subsection{\@tocline{2}{0pt}{2.3em}{}{}}
\makeatother

\usepackage{amsmath,amsfonts,amssymb,amsthm,amscd}

\usepackage{epsfig}
\usepackage{psfrag}
\usepackage{perpage}
\usepackage{url}
\usepackage{color}
\usepackage[english]{babel}
\usepackage{bbm}
\usepackage{epstopdf}
\usepackage[utf8]{inputenc}
\usepackage[T1]{fontenc}
\usepackage{microtype}
\usepackage{hyperref}
\hypersetup{colorlinks=true, linkcolor=blue, citecolor=blue,  filecolor=blue, urlcolor=blue}
\usepackage{verbatim}
\usepackage{subcaption}
\usepackage[ruled,vlined]{algorithm2e}
\usepackage{multirow}
\usepackage{booktabs}

\counterwithin{figure}{section}
\numberwithin{figure}{section}

\numberwithin{equation}{section}

\usepackage{fullpage}

\newcommand{\ee}{\mathbb{E}}

\newcommand{\pp}{\mathbb{P}}

\renewcommand{\hat}{\widehat}

\newcommand{{\paa}[1]}{p_{00,#1}}
\newcommand{{\pab}[1]}{p_{01,#1}}
\newcommand{{\pba}[1]}{p_{10,#1}}
\newcommand{{\pbb}[1]}{p_{11,#1}}

\begin{document}


\title{Parameter estimation in interacting particle systems \\ on dynamic random networks}
\thanks{The work of JW was supported by the European Union’s Horizon 2020 research and innovation programme under the Marie Skłodowska-Curie grant agreement no.\ 101034253; SB received the support from the same grant while affiliated with Leiden University. SB was further supported through “Gruppo Nazionale per l’Analisi Matematica, la Probabilità e le loro Applicazioni” (GNAMPA-INdAM). The authors are grateful to Michel Mandjes for useful and fruitful discussions.
  \includegraphics[height=1em]{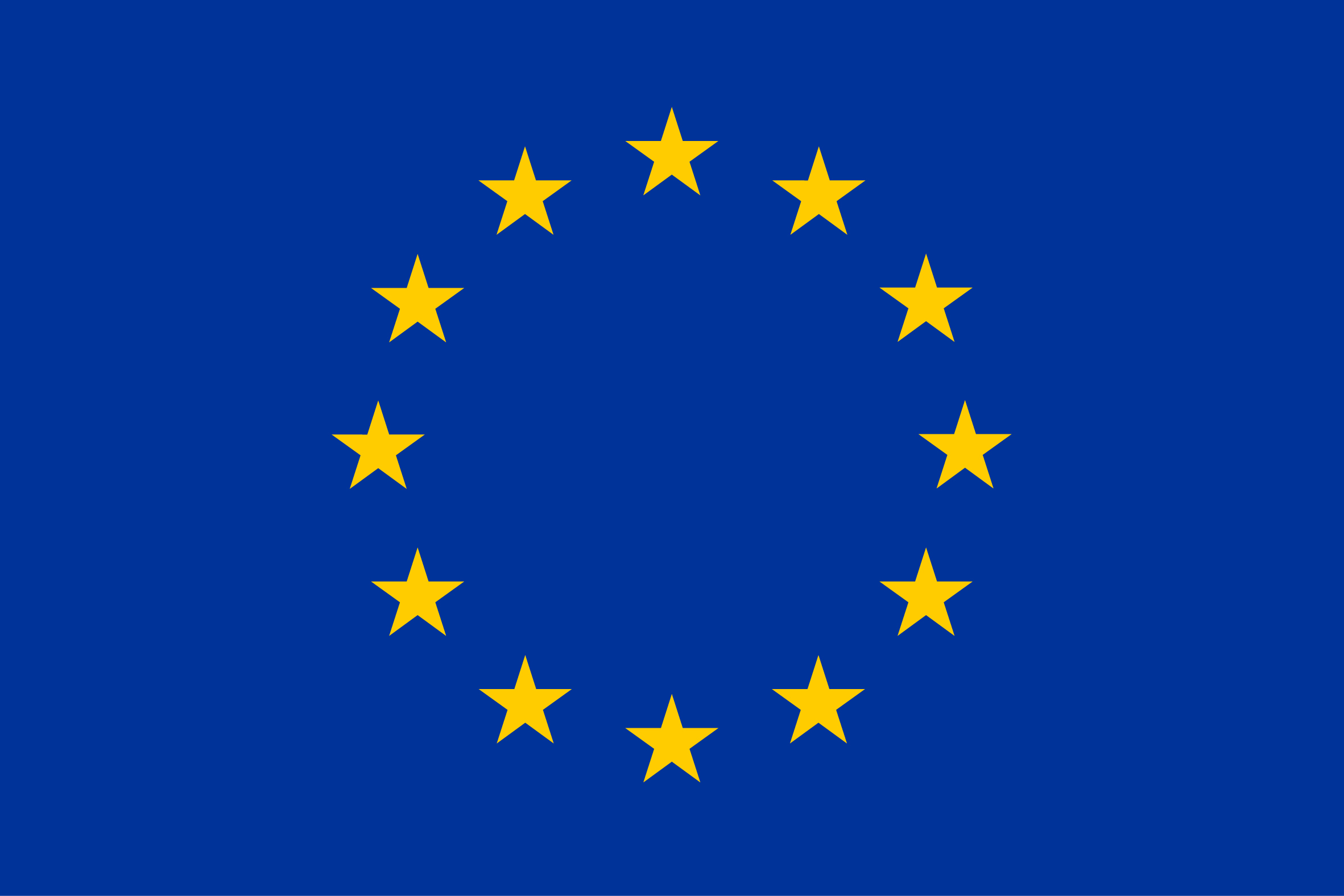}}%

\author{Simone Baldassarri}
\altaffiliation{Gran Sasso Science Institute, Viale Francesco Crispi 7, 67100 L’Aquila, Italy}
\email{simone.baldassarri@gssi.it} 

\author{Jiesen Wang}
\altaffiliation{Korteweg-de Vries Institute for Mathematics, University of Amsterdam, Science Park 904, 1098 XH Amsterdam, The Netherlands}
\email{j.wang2@uva.nl}


\begin{abstract}
In this paper we consider a class of interacting particle systems on dynamic random networks, in which the joint dynamics of vertices and edges acts as {\it one-way feedback}, i.e., edges appear and disappear over time depending on the state of the two connected vertices, while the vertex dynamics does not depend on the edge process. Our goal is to estimate the underlying dynamics from partial information of the process, specifically from snapshots of the total number of edges present. We showcase the effectiveness of our inference method through various numerical results.
\end{abstract}

\maketitle



\section{Introduction}

In recent years, the study of interacting particle systems (IPS) on {\it dynamic} random networks has become a central framework for understanding the evolution of complex systems. Classical IPS models, such as the voter model, contact process, and exclusion process, typically assume a {\it static} interaction network. However, this assumption often fails to reflect the temporal and adaptive structure of real--world systems, where both the states of individual units and the network of interactions evolve concurrently. This co--evolutionary dependence gives rise to feedback loops and emergent behavior that challenge traditional modeling and inference techniques.

To address this, we consider IPS on dynamic random graphs, where edges appear and disappear over time, potentially depending on the internal states of the nodes. We focus on a class of models with {\it unidirectional dependence}: the states of the nodes influence the dynamics of the edges, but not viceversa. This structure offers a simplification that remains both realistic and analytically tractable, particularly in systems where node behavior is shaped by external inputs or intrinsic rules. Our central contribution is to develop a parametric inference framework for such systems under minimal observability assumptions. Specifically, we assume that only the total number of active edges is observed at each discrete time point.

We use the method of moments (see e.g.\ \cite{HMW25,MW24,MW25} for recent applications to dynamic random graphs), which avoids the intractability of likelihood methods in partially observed and non--Markovian settings. By leveraging moment equations that relate expected edge counts and their covariances to model parameters, we construct estimators that are both easy to compute and statistically meaningful. Numerical experiments suggest that these estimators are consistent and asymptotically normal, supporting their use for inference in complex evolving networks. This general framework is applicable to a wide range of domains where the interaction structure is shaped by local node--level dynamics. We mention below three main applications of our model.

\begin{itemize}
\item {\it Opinion Dynamics.} In the context of voter models, individuals hold discrete opinions that evolve via interactions with their social contacts. While classical models assume a fixed network, empirical studies show that social ties often change in response to opinion alignment. Individuals may distance themselves from those with opposing views or become more connected to like--minded peers, introducing a feedback loop between opinions and network structure. 
This interplay has been extensively explored in the literature through numerical experiments (see e.g.\ \cite{AVHG25,DGLMSSV12,GZ06,HN06,LHAJS20,PTN06,PMDDC24,RG17,ZG06}), while only recently through theoretical analyses (see e.g.\ \cite{AdHR25,ABHdHQ25,BBdHM24,BK25,BS17}), under frameworks incorporating either one--way or two--way feedback between opinion dynamics and network evolution. In such systems, one often observes only aggregate opinion fractions or local counts of active interactions rather than full microstate trajectories. This partial knowledge makes it natural to estimate interaction parameters under partial information rather than from complete history of the network.

\item{\it Communication Networks.} In systems such as the Internet of Things (IoT), peer--to--peer architectures, or mobile ad hoc networks (MANETs), nodes represent devices with internal activity states (transmitting or idle), and communication links are created or removed depending on these states. Similar principles underlie real--world protocols such as CSMA/CA (see e.g.\ \cite{BBM06,DeGRACD16,DDT09,LDPF22,ZAWR25}) and DTNs (see e.g.\ \cite{MKMV23,WJCL23,ZAZ04}), where access and connectivity depend on local activity and mobility.  
Complementary investigations in dynamic and temporal networks highlight how time--varying connectivity, local sensing, and intermittent interactions affect reachability, information spreading, and transmission delays (see e.g.\ \cite{FFSR24,JRDRK17,TCS17}). Together, these works reveal how device activity and connectivity shape performance in dynamic networks. Measurements typically record only the presence or activity of connections at discrete times, rather than the full temporal evolution of the network. Consequently, parameter inference must rely on these snapshot--based moments of active edges, highlighting the relevance of our framework under partial observability.

\item
{\it Neural Systems.} Brain networks exhibit complex feedback between neural activation and structural connectivity. Synapses may strengthen or weaken based on recent activity (as in Hebbian plasticity), while the network topology constrains possible activation patterns. Empirical studies have shown how task--related neural activity reshapes functional connectivity and center organization, highlighting the dynamic nature of the structure of the brain network (see, e.g. \cite{TLVMGB16,YDR23,YI23}). On the other hand, theoretical work has demonstrated how spike--timing--dependent plasticity can give rise to emergent network topologies, such as clustering and modularity, even from initially random networks (see e.g.\ \cite{BBLBICAB17, RV21, SLP19, TB16,vKLdAH18, YWDLJ24}). Our framework applies to models where neurons are nodes with dynamic states, and synaptic connections evolve in response to co--activation history. In this setting, data are often available only through aggregated observables, making inference under partial observability both necessary and challenging. 
\end{itemize}

To the best of our knowledge, this paper is the first that statistically infers interacting particle systems on dynamic random graphs with one--way feedback under partial information of the process. Despite its analytical tractability, the model is highly flexible and can be extended to capture more intricate mechanisms, including protocols with mutual dependencies between the dynamics of vertices and edges. As such, it provides a foundational step toward the analysis of more complex co--evolutionary systems where vertex states and network structure interact through two--way feedback. This framework opens broad opportunities for data--driven investigation across social, technological, and biological domains, especially in settings where only aggregate or incomplete data are available.

The rest of the paper is organized as follows. In Section \ref{sec:model}, we introduce the model of our interest and we discuss how it relates to opinion dynamics, communication networks and neural systems. In Section \ref{sec:methodology}, we derive the moment equations that we need in order to estimate the underlying dynamics, while Section \ref{sec:simulations} provides numerical experiments showing that our estimators are consistent and asymptotically normal. In Section \ref{sec:conclusions}, we draw our conclusions. Finally, Appendices \ref{sec:appA}, \ref{sec:appB} and \ref{sec:appC} contain some technicalities that arise in Section \ref{sec:methodology}.

\section{Definition of the model}
\label{sec:model}

We consider a model with a one--way interaction, i.e., the states of the vertices are independent of the states of the edges, but the reverse is not true. At any time $t\geq0$, the set of vertices is $V=[n]=\{1,2,...,n\}$ and any vertex holds state either $+$ or $-$. We assume that each edge is initially active with probability $p_0$, independently of everything else (i.e., we start from an Erdös--Rényi random graph on $n$ vertices with edge retention probability $p_0$), and half vertices hold state $+$, while the others hold state $-$.
Let $E(t)$ be the set of active edges at time $t$ and 
\[
x_i(t) = \mbox{state of vertex $i$ at time $t$}.
\]
At each discrete time step, given $\alpha\in(0,1)$, the model is characterised by the following dynamics:
\begin{itemize}
\item with probability $\alpha$, a vertex is selected uniformly at random and it changes state;
\item with probability $1-\alpha$, each edge $ij$, with $i<j$, is active with probability $\pi_+$ (resp.\ $\pi_-$) if $x_i(t)=x_j(t)=+$ (resp.\ $x_i(t)=x_j(t)=-$), and with probability $f(\pi_+,\pi_-)$ if $x_i(t)$ and $x_j(t)$ are different, where $\pi_+,\pi_-\in[0,1]$ and $f:[0,1]^2\to[0,1]$ is a deterministic prescribed function.
\end{itemize}
In applications, the parameter $\alpha$ should be in general thought of a very low value, because vertices only rarely change their state randomly.
By definition, it is immediate to see that the vertex dynamics is independent of the edges dynamics, but the reverse is not true. 
The method we present in this paper does not depend on the particular choice of the function $f$, which can be different depending on the system we aim to model. Below, we explain how this dynamics is interpreted in our three illustrative applications, what are the features the function $f$ should have in each of these examples, and how, in real--world settings, information about the network’s connectivity is typically available through observable interactions.

\medskip
{\it Opinion Dynamics.} Vertices represent individuals or agents in a social system, and their state is seen as {\it supportive}/{\it opposing}, i.e., binary stances on a given issue. Specifically, vertices tend to align and interact primarily with others in the same state, while interactions with individuals holding the opposing view are less likely. Thus, we require that the function $f$ satisfies $f(\pi_+,\pi_-)<\min\{\pi_+,\pi_-\}$. A possible choice is then
\begin{equation}\label{eq:fdef1}
f(\pi_+,\pi_-)=
\begin{cases}
    \dfrac{\pi_+\pi_-}{\pi_++\pi_-} &\text{ if } \pi_++\pi_-\neq0, \\
    0 &\text{ if } \pi_+=\pi_-=0,
\end{cases}
\end{equation}
which attains the maximum $1/2$ when $\pi_+=\pi_-=1$. In this context, edges represent social ties or channels of influence (e.g., conversations, online interactions). Active ties may be inferred from observed patterns of communication, or interaction frequency among individuals, which indicate potential pathways for influence. The dynamics can then be interpreted as follows. A vertex update occurs when an individual changes opinion in response to external factors such as media, personal beliefs, or intrinsic tendencies, reflecting the cumulative impact of influences shaping decisions over time, while edge updates represent the adaptation of social connections, that is, the tendency to form or break ties based on shared opinions. The probability of vertex transition $\alpha$ plays the role of {\it social susceptibility}, which quantifies the response of the average opinion (magnetization) to an external bias (e.g., $\alpha$ increases with social reinforcement or persuasive messaging).

\medskip
{\it Communication Networks.} Vertices represent devices or users in a communication network, and their state is seen as {\it transmitting}/{\it idle}, i.e., binary communication activity modes. Specifically, a vertex in the transmitting (resp.\ idle) state is actively sending data (resp.\ silent) and is more (resp.\ less) likely to engage in or maintain communication links. For example, two transmitting devices are more likely to exchange messages than a pair where one or both are idle, and constraints like energy or mobility affect connections. Thus, we require that the function $f$ satisfies $\pi_-<f(\pi_+,\pi_-)<\pi_+$. A possible choice is then
\begin{equation}\label{eq:fdef2}
f(\pi_+,\pi_-)=\dfrac{\pi_++\pi_-}{2}, \quad \pi_+>\pi_-.
\end{equation}
In this context, edges represent communication channels or data links (e.g., packet routing paths). Active connections are typically captured through logs of message exchanges, connection durations, or data traffic volumes between devices, providing a proxy for the underlying network activity. The dynamics can be then interpreted as follows. A vertex update occurs when a device starts or stops transmitting due to traffic demand or external interference, while edge updates represent the reconfiguration of communication links, such as link establishment or failure recovery. The probability of vertex transition $\alpha$ plays the role of {\it transmission likelihood} (e.g., $\alpha$ increases with data load or priority level).

\medskip
{\it Neural Systems.} 
Vertices represent neurons, whose states are modeled as {\it active} or {\it quiescent}, corresponding to binary firing modes. The BOLD (Blood Oxygenation Level Dependent) signal is the primary contrast mechanism used in functional Magnetic Resonance Imaging (fMRI) to infer brain connectivity. However, BOLD is an indirect measure of neural activity, relying on changes in blood oxygenation; active neurons tend to produce a stronger BOLD signal~\cite{LPATO01}.
A good choice for the function $f$ is then given, for instance, by \eqref{eq:fdef2}.
In this context, edges represent functional connectivity derived from fMRI data. The dynamics can then be interpreted as follows: a vertex update occurs when a neuron flips its state in response to external input or noise, while edge updates reflect changes in functional connectivity. The vertex transition probability~$\alpha$ characterizes {\it neuronal excitability} (e.g., $\alpha$ may increase under the influence of neuromodulators such as dopamine).

\medskip
In all the three cases, the active connections in the network reflect the dynamic substrate over which information or influence flows, and they are often the most accessible features in empirical data. 

To illustrate the behavior of the model, we now focus on the two following representative parameter choices:
\begin{itemize}
    \item[(i)] $\pi_+=\pi_-$, where we refer to the common value as $\pi$;
    \item[(ii)] $|\pi_+-\pi_-|\approx 1$, i.e., one of the two parameters is very close to 1 and the other to 0.
\end{itemize}
In (i), the probability of having an active edge between vertices having the same state is equal, while for vertices having different states this probability depends on the choice of $f$. In particular, $f$ in \eqref{eq:fdef1} and \eqref{eq:fdef2} reduces to $\pi/2$ and $\pi$, respectively, so that in the latter case, at each time step in which the edges are resampled, the network is an Erdös-Rényi random graph on $n$ vertices with edge retention probability $\pi$.

In (ii), assume that $\max\{\pi_+,\pi_-\}=\pi_+$ without loss of generality. Thus, there is a high (resp.\ low) probability of having an active edge between vertices having both same state $+$ (resp.\ $-$). For vertices having different states, $f$ in \eqref{eq:fdef1} and \eqref{eq:fdef2} is close to $0$ and $1/2$, respectively.

In \eqref{eq:fdef1}--\eqref{eq:fdef2} we choose $f$ as a symmetric function in the parameters $\pi_+$ and $\pi_-$. This is a natural choice when the rule for forming an edge depends only on the states’ frequencies and treats the two states equally. A typical modeling motivation is when there is no a priori bias (i.e., no leadership, no external field). However, when vertices with different states have different behavioral roles, an asymmetric function becomes meaningful. In fact, in social networks some people (say $+$) might be more likely to attract other (say $-$) people (e.g., during a marketing campaign). In information--exchange systems, asymmetric edge probabilities may appear because sending and receiving often require different amounts of effort or resources. Finally, in neuroscience, active and quiescent neurons have different connectivity propensities, so that edge--formation should weight the state that is more attractive more heavily. 

A natural way to generate the asymmetric version of $f$ as in \eqref{eq:fdef1} is to think that one of the two states has to do {\it more work} than the other to make the connection possible. For instance, assume that the edge creation requires $k_+$ (resp.\ $k_-$) units of work from a $+$ (resp.\ $-$) vertex, then the effected edge--formation probability is given by
\begin{equation}\label{eq:fdef1asy}
f(\pi_+,\pi_-) = \dfrac{1}{\frac{k_+}{\pi_+}+\frac{k_-}{\pi_-}} = \dfrac{\pi_+\pi_-}{k_+\pi_-+k_-\pi_+}.
\end{equation}
As an example, in opinion dynamics a $+$ vertex may require multiple exposures to a $-$ neighbor before it is willing to form a connection (e.g., to listen), but the $-$ vertex may require only one.

A convenient way to generate the asymmetric version of $f$ as in \eqref{eq:fdef2} is to imagine weighted sampling of the two endpoints. Suppose edge--formation between a randomly chosen pair occurs by:
\begin{itemize}
    \item[(i)] A vertex of state $s$ contributes a {\it strength} $\beta_s$ toward forming an edge when involved, with $\beta_s=a_s\pi_s$, where $a_s>0$ is a {\it weight} measuring how strongly that parameter is expressed (activity or number of connection attempts).
    \item[(ii)] When two vertices with different states meet, the pairwise formation propensity is the sum of the contributions from the two endpoints, i.e., it equals $\beta_++\beta_-=a_+\pi_+ + a_-\pi_-$.
\end{itemize}
The attachment probability is then given by
\begin{equation}\label{eq:fdef2asy}
f(\pi_+,\pi_-) = \dfrac{a_+\pi_+ + a_-\pi_-}{a_++a_-},
\end{equation}
where the case $a_+=a_-$ gives the attachment function in \eqref{eq:fdef2}, otherwise one of the two state has stronger attractiveness.
This can be understood as follows. In communication networks, the connection probability is driven more strongly by how active the transmitter is, with the idle node playing a secondary role. In neural systems, to establish a synapse, an active neuron might require multiple spikes (work units) to release enough neurotransmitter, while a quiescent neuron needs only one spike.

\section{Methodology}
\label{sec:methodology}
Let 
\[
N(t) = \hbox{number of vertices with state + at time } t,
\]
and
\begin{equation}\label{eq:defS}
\begin{array}{ll}
S(t) &= \hbox{number of active edges at time } t \\
&= \displaystyle\sum_{i<j} \mathbf{1}\{(i,j)\in E(t)\}.
\end{array}
\end{equation}

We assume that only the number of active edges, denoted by the sequence $\{S(t)\}_{t = 1, \ldots, K}$, is observable, where $K$ is the number of observations. Our objective is to estimate the underlying parameter set $\{ \pi_+, \, \pi_-, \, \alpha \}$ that governs the network dynamics, using the observed history of edge counts over time. A central challenge in this setting is the limited observability: while we have access to the total number of edges at each time step, we do not observe the identities of individual edges.

When the total number of active edges changes between two time steps, we can infer with certainty that an edge update has occurred. However, if the edge count remains unchanged, the situation is more ambiguous. The update may involve a vertex whose rewiring preserves the total number of edges, or it may be an edge update that results in the removal and addition of the same number of edges. Consequently, the observed edge count alone does not uniquely determine the type of update that has taken place.

To address this challenge and extract information from the available data, we adopt a method--of--moments approach. By linking the expected value and higher--order moments of the edge count process to the model parameters, we derive estimators that can be matched to empirical moments computed from the observations. A key consideration in this framework is the choice of which moments to use. The selected moments should satisfy the following criteria:
\begin{itemize}
    \item An expression (not necessarily in closed form) for each chosen moment must be available under the model.
    \item The empirical counterparts of the selected moments must be computable from the observed data.
    \item The number of linearly independent moment equations should match the number of unknown parameters.
\end{itemize}

In what follows, we provide the three moment equations we will use for the estimate of the parameters $\pi_+$, $\pi_-$ and $\alpha$.

We assume that our system starts at {\it stationarity}, i.e., the distributions of $N(t)$ and $S(t)$ do not depend on $t\in[0,K]$. Thus, when considering
\[
P_{t,k} = \mathbb{P} (N(t) = k),
\]
we may remove the time dependence. By the definition of the double dynamics, we know that, in order that the system reaches $k$ vertices with state $+$ in a single transition, with $1\leq k \leq n-1$, there are three possibilities:
\begin{itemize}
    \item[(i)] the system has $k-1$ vertices with state $+$, the vertex dynamics is selected and a vertex changes its state from $-$ to $+$, which occurs with probability $\alpha \frac{n-k+1}{n}P_{k-1}$;
    \item[(ii)] the system has $k+1$ vertices with state $+$, the vertex dynamics is selected and a vertex changes its state from $+$ to $-$, which occurs with probability $\alpha \frac{k+1}{n} P_{k+1}$;
    \item[(iii)] the system has $k$ vertices with state $+$ and the edge dynamics is selected, which occurs with probability $(1-\alpha) P_{k}$.
\end{itemize}
Thus, at stationarity, the probability distribution of the vertex process obeys the balance equation (stationary master equation) \cite{G2009,vK2007}, which, for $1\leq k \leq n-1$, is given by
\[
P_{k} = \alpha \dfrac{n-k+1}{n}P_{k-1} + \alpha \dfrac{k+1}{n} P_{k+1} + (1-\alpha) P_{k},
\]
thereby leading to
\begin{equation}\label{eq:recequation}
    P_k = \dfrac{n-k+1}{n}P_{k-1} + \dfrac{k+1}{n} P_{k+1},
\end{equation}
with boundary conditions $P_1=nP_0$ and $P_{n-1}=nP_n$. Indeed, $P_0$ and $P_n$ satisfy the equations
\[
P_0 = \dfrac{\alpha}{n} P_1 + (1 - \alpha) P_0, \quad P_n = \dfrac{\alpha}{n} P_{n-1} + (1 - \alpha) P_n.
\]
In what follows, we will need the explicit characterization of $P_k$, which is given by
    \begin{equation}\label{eq:recurrence}
        P_k= \dfrac{1}{2^n} \binom{n}{k}, \qquad 0\leq k\leq n.
    \end{equation}
The proof of \eqref{eq:recurrence} is deferred to Appendix \ref{sec:appA}. 
Since
\[
\ee[S(t)] = \displaystyle\sum_{i<j} \pp\left( (i,j)\in E(t) \right),
\]
we deduce that, for any $t\geq1$,
\[
\begin{array}{ll}
\ee[S(t)] - \alpha\ee[S(t-1)] \\
\quad = (1-\alpha) \displaystyle\sum_{i<j} \sum_{k=0}^n \pp\left( (i,j)\in E(t) | N_+(t)=k \right) P_k \\
\quad = (1-\alpha)\displaystyle\sum_{k=0}^n F(k,\pi_+,\pi_-) P_k,
\end{array}
\]
where
\[
\begin{array}{ll}
F(k,\pi_+,\pi_-) \\
\quad\quad = \displaystyle\binom{k}{2}\pi_+ + \binom{n-k}{2}\pi_- + k(n-k) f(\pi_+,\pi_-).
\end{array}
\]
This, at stationarity, gives
\begin{equation}\label{eq:moment1}
\ee[S(t)]= \sum_{k=0}^n F(k,\pi_+,\pi_-) P_k.
\end{equation}

In what follows, we can write the explicit characterization of $\ee[S^2(t)]$, with $t\in[0,K]$, which is given by
    \begin{equation}\label{eq:secondmoment}
    \ee[S^2(t)]=\ee[S(t)]+\sum_{k=0}^n (E_k^1 + E_k^2 + E_k^3) P_k,
    \end{equation}
    where 
    \begin{equation*}
    \renewcommand{\arraystretch}{2}
    \begin{array}{ll}
    E_k^1 &= \displaystyle\binom{k}{2}\pi_+ \Bigg\{ \left[ \binom{k}{2} - 1 \right] \pi_+ + \binom{n-k}{2}\pi_- \\
    &\quad \displaystyle + k (n-k) f(\pi_+,\pi_-)\Bigg\}, \\
    E_k^2 &= \displaystyle\binom{n-k}{2}\pi_-  \Bigg\{ \binom{k}{2} \pi_+ + \left[\binom{n-k}{2}-1\right]\pi_- \\
    &\quad \displaystyle + k (n-k) f(\pi_+,\pi_-)\Bigg\}, \\
    E_k^3 &= k(n-k) f(\pi_+,\pi_-) \displaystyle \Bigg \{ \binom{k}{2} \pi_+ + \binom{n-k}{2}\pi_- \\
    &\quad \displaystyle + \left[k (n-k) - 1 \right] f(\pi_+,\pi_-)\Bigg\}.
    \end{array}
    \end{equation*}

The proof of \eqref{eq:secondmoment} is deferred to Appendix \ref{sec:appB}.

Observe that the expressions for $\mathbb{E}[S(t)]$ and $\mathbb{E}[S^2(t)]$ depend only on the parameters $\pi_+$ and $\pi_-$, but not on $\alpha$. Moreover, $\mathbb{E}[S(t)]$ and $\mathbb{E}[S^2(t)]$ are linearly independent, in the sense that one cannot be expressed as a function of the other. Also, the empirical counter parts can be expressed as 
\[
M_{1,k} := \frac{1}{K} \sum_{t = 1}^{K} S(t), \qquad M_{2,k} := \frac{1}{K} \sum_{t = 1}^{K} S^2(t) \,,
\]
respectively.
This means that $\pi_+$ and $\pi_-$ can be separately estimated from $\alpha$ by solving
\begin{equation} \label{eq:estimation1}
    \mathbb{E}[S(t)] = M_{1,k}, \qquad \mathbb{E}[S^2(t)] = M_{2,k} \,.
\end{equation}
The advantage of separate estimation is the simplicity of the procedure, which makes the estimation process more efficient in terms of computational speed (see \cite{MW25}). We stress that, if the link function $f$ is symmetric in the two variables and $(\pi_+,\pi_-)=(\pi_1,\pi_2)$ is a solution of the moment equations \eqref{eq:estimation1}, then due to the nature of the model also $(\pi_+,\pi_-)=(\pi_2,\pi_1)$ is. We prove this in Appendix \ref{sec:appC}. This does not happen if the function $f$ is asymmetric, as discussed in Section \ref{sec:simulations} below.

Finally, to estimate $\alpha$ we require a moment that depends on $\alpha$, which involves tracking the evolution of the system over time. To this end, we proceed as follows. Consider the joint process of number of vertices having state $+$ and the number of active edges at time $t$, namely the process $(N(t),S(t))$. Note that this is a Markov chain, but at each time step we observe only $S(t)$, which however is not a Markov chain. The idea is then the following. 

First, write down the transition probabilities $P_{(i,k)(j,\ell)}$ for $i,j=0,...,n$ and $k,\ell=0,...,\binom{n}{2}$, so that it is possible to compute the (unique) stationary distribution $\{\pi(i,k) : i=0,...,n \text{ and } k=0,...,\binom{n}{2}\}$. The transition probabilities are defined as follows:
\begin{equation}
     P_{(i,k)(j,\ell)}= 
     \begin{cases}
         0 &\text{ if } |j-i|\geq2, \\
         0 &\text{ if } |j-i|=1, k\neq\ell, \\
          \alpha(n-i)/n &\text{ if } j=i+1, i\neq n, k=\ell, \\
          \alpha i/n &\text{ if } j=i-1, i\neq0, k=\ell, \\
         (1-\alpha) P^*_{i\ell} &\text{ if } j=i,
     \end{cases}
\end{equation}
where
\[
\begin{array}{ll}
P^*_{i\ell} = \\
 \displaystyle\sum_{r_1=0}^{r_1^*} \sum_{r_2=0}^{r_2^*} \binom{\binom{i}{2}}{r_1} \binom{\binom{n-i}{2}}{r_2} \binom{i(n-i)}{\ell-r_1-r_2} \\
 \times \pi_+^{r_1} (1-\pi_+)^{\binom{i}{2}-r_1} \pi_-^{r_2} (1-\pi_-)^{\binom{n-i}{2}-r_2} \\
 \times \left(f(\pi_+,\pi_-)\right)^{\ell-r_1-r_2} \left(1-f(\pi_+,\pi_-)\right)^{i(n-i)-\ell+r_1+r_2},
\end{array}
\]
with $r_1^*=\min\left\{\ell, \binom{i}{2}\right\}$ and $r_2^*=\min\left\{\ell-r_1, \binom{n-i}{2}\right\}$. Note that the term $P^*_{i\ell}$ corresponds to the probability that, once the edge dynamics is selected, the resulting model has $\ell$ active edges, where the sums over $r_1$ and $r_2$ correspond to the events where the number of active edges connecting vertices having both states $+$ and $-$ are $r_1$ and $r_2$, respectively. As an example, when $n=3$, $\alpha=0.3$, $\pi_+=0.9$, $\pi_-=0.4$ and $f$ as in \eqref{eq:fdef2}, the transition matrix is given in Figure \ref{fig:P}.
\begin{figure*}
\centering
\[
{\large
P=\left(
\begin{smallmatrix}
    0.1512 & 0.3024 & 0.2016 & 0.0448 & 0.3 & 0 & 0 & 0 & 0 & 0 & 0 & 0 & 0 & 0 & 0 & 0 \\
    0.1512 & 0.3024 & 0.2016 & 0.0448 & 0 & 0.3 & 0 & 0 & 0 & 0 & 0 & 0 & 0 & 0 & 0 & 0 \\
    0.1512 & 0.3024 & 0.2016 & 0.0448 & 0 & 0 & 0.3 & 0 & 0 & 0 & 0 & 0 & 0 & 0 & 0 & 0 \\
    0.1512 & 0.3024 & 0.2016 & 0.0448 & 0 & 0 & 0 & 0.3 & 0 & 0 & 0 & 0 & 0 & 0 & 0 & 0 \\
    0.1 & 0 & 0 & 0 & 0.0514 & 0.2254 & 0.3049 & 0.1183 & 0.2 & 0 & 0 & 0 & 0 & 0 & 0 & 0 \\
    0 & 0.1 & 0 & 0 & 0.0514 & 0.2254 & 0.3049 & 0.1183 & 0 & 0.2 & 0 & 0 & 0 & 0 & 0 & 0 \\
    0 & 0 & 0.1 & 0 & 0.0514 & 0.2254 & 0.3049 & 0.1183 & 0 & 0 & 0.2 & 0 & 0 & 0 & 0 & 0 \\
    0 & 0 & 0 & 0.1 & 0.0514 & 0.2254 & 0.3049 & 0.1183 & 0 & 0 & 0 & 0.2 & 0 & 0 & 0 & 0 \\
    0 & 0 & 0 & 0 & 0.2 & 0 & 0 & 0 & 0.0086 & 0.1090 & 0.3162 & 0.2662 & 0.1 & 0 & 0 & 0 \\
    0 & 0 & 0 & 0 & 0 & 0.2 & 0 & 0 & 0.0086 & 0.1090 & 0.3162 & 0.2662 & 0 & 0.1 & 0 & 0 \\
    0 & 0 & 0 & 0 & 0 & 0 & 0.2 & 0 & 0.0086 & 0.1090 & 0.3162 & 0.2662 & 0 & 0 & 0.1 & 0 \\
    0 & 0 & 0 & 0 & 0 & 0 & 0 & 0.2 & 0.0086 & 0.1090 & 0.3162 & 0.2662 & 0 & 0 & 0 & 0.1 \\
    0 & 0 & 0 & 0 & 0 & 0 & 0 & 0 & 0.3 & 0 & 0 & 0 & 0.0007 & 0.0189 & 0.1701 & 0.5103 \\
    0 & 0 & 0 & 0 & 0 & 0 & 0 & 0 & 0 & 0.3 & 0 & 0 & 0.0007 & 0.0189 & 0.1701 & 0.5103 \\
    0 & 0 & 0 & 0 & 0 & 0 & 0 & 0 & 0 & 0 & 0.3 & 0 & 0.0007 & 0.0189 & 0.1701 & 0.5103 \\
    0 & 0 & 0 & 0 & 0 & 0 & 0 & 0 & 0 & 0 & 0 & 0.3 & 0.0007 & 0.0189 & 0.1701 & 0.5103
\end{smallmatrix}
\right),
}
\]
\caption{Example of transition matrix when $n=3$, $\alpha=0.3$, $\pi_+=0.9$, $\pi_-=0.4$ and $f$ as in \eqref{eq:fdef2}.}
\label{fig:P}
\end{figure*}
where the states are ordered lexicographically. After computing the stationary distribution, it is possible to compute all the cross--moments as
\[
\begin{array}{ll}
\ee\left[S(t)S(t+k)\right] \\
\quad= \displaystyle\sum_{\ell,r} \ell r \mathbb{P}(S(t)=\ell, S(t+k)=r) \\
\quad= \displaystyle\sum_{\ell,r} \ell r \sum_{i}\mathbb{P}(N(t)=i,S(t)=\ell, S(t+k)=r) \\
\quad= \displaystyle\sum_{i,j,\ell,r} \ell r P^k_{(i,\ell)(j,r)} \pi(i,\ell), 
\end{array}
\]
where $k\geq1$ and the last equality holds for any $t$ because our system starts at stationarity. In particular, for $k=1$ we get
\[
\ee\left[S(t)S(t+1)\right] = \sum_{i,j,\ell,r} \ell r P_{(i,\ell) (j,r)} \pi(i,\ell),
\]
thereby allowing us to write the one--lag mean squared difference
\[
\begin{array}{ll}
\ee[(S(t+1)-S(t))^2] \\
\quad = \ee[S(t+1)^2] + \ee[S(t)^2] - 2 \ee[S(t)S(t+1)] \\
\quad = 2 \ee[S(t)^2] - 2 \displaystyle\sum_{i,j,\ell,r} \ell r P_{(i,\ell) (j,r)} \pi(i,\ell).
\end{array}
\]
To estimate the parameter $\alpha$, we will use the above expression for $\ee[(S(t+1)-S(t))^2]$. Note that the last quantity can be written in terms of $\text{Cov}(S(t),S(t+1))$ as
\[
\begin{array}{ll}
\ee[(S(t+1)-S(t))^2] \\
\qquad \quad \quad= 2\text{Var}(S(t)) - 2\text{Cov}(S(t),S(t+1)).
\end{array}
\]
Below we give a heuristic argument to explain the reason why we estimate the parameter $\alpha$ by using $\ee[(S(t+1)-S(t))^2]$ and not simply $\ee\left[S(t)S(t+1)\right]$. The cross--moment $\ee\left[S(t)S(t+1)\right]$ captures temporal correlation between the number of active edges at time $t$ and at time $t+1$. The cross--moment can be relatively {\it insensitive} in $\alpha$ because, even when edge resampling occurs, the vertex states might not change enough to affect many edge activations. Similarly, vertex updates may not shift global structure drastically in one time step. Thus, different values of $\alpha$ can yield similar values of the cross--moment, reducing its usefulness for parameter estimation. When $\alpha$ approaches 0, the system is dominated by edge resampling: the vertex states stay almost static over time, and at each time step, the entire edge set is redrawn independently of the previous configuration (but still based on the same state snapshot). Thus, $S(t)$ and $S(t+1)$ are both random variables drawn from similar distributions, but conditionally independent given the vertex states. Hence,
\[
\ee[S(t)S(t+1)] = \ee[S(t)]^2 + (\text{small correction term)}.
\]
This leads to a value close to the square of the mean, and does not change much with small variations in $\alpha$, because the dominant mechanism (edge resampling) remains. So, the moment is almost flat as $\alpha$ varies near zero. When $\alpha$ approaches 1, the system is dominated by vertex resampling. Since vertex states evolve slowly (only one vertex flips at each time step), from $t$ to $t+1$ the vertex configuration changes only slightly and the change in active edges is minor. Thus, $S(t+1)\approx S(t)$, and the correlation becomes $\ee[S(t)S(t+1)]\approx\ee[S(t)^2]$. Again, small changes in $\alpha$ around 1 have little effect, because the number of vertex states that change per step remains close to 1. The moment then stays flat as $\alpha$ varies near one. For intermediate $\alpha$, both the vertex and edge dynamics are present. Even in this regime, the cross--moment aggregates across many possible outcomes. Because edge activity is a combination of slowly changing vertex states and sudden edge resets, the resulting $S(t+1)$ can still be mildly correlated with $S(t)$ depending on how much vertex state shift occurred recently. Hence, the dependence of $\ee[S(t)S(t+1)]$ on $\alpha$ tends to be relatively flat over a large central region. In practice, the derivative of this moment with respect to $\alpha$ is small across much of the domain. This makes it a bad candidate for parameter estimation, because different values of $\alpha$ yield similar values of $\ee[S(t)S(t+1)]$ and this reduces identifiability. On the other hand, at stationarity, the moment $\ee[(S(t+1)-S(t))^2]$ represents how variable the system is from one time step to the next. With similar arguments as above, we can see that when $\alpha$ approaches zero, the variance of $S(t+1)-S(t)$ reflects only edge-level randomness (binomial fluctuations), so that $\ee[(S(t+1)-S(t))^2]$  captures only edge fluctuation noise. When $\alpha$ approaches one, $S(t)$ can change substantially when a vertex flips and causes many edges to switch state (depending on its neighbors’ states). In this regime, the evolution of the vertex states drives non--trivial changes in $S(t)$, and the squared increment captures this. For intermediate $\alpha$, the net change in $S(t)$ reflects both how stable vertex states are, and how much randomness the edge layer introduces. Thus, $\ee[(S(t+1)-S(t))^2]$ blends both sources of fluctuation and is more sensitive to their balance, which is governed by $\alpha$. All this implies that, as $\alpha$ varies, the dominant source of variability shifts, thereby causing a non--flat, more sensitive response in the moment.

The covariance can also be used to estimate $\alpha$, as it reflects the correlation between consecutive steps and partially captures edge fluctuations. However, compared to the expression $\mathbb{E}[(S(t+1) - S(t))^2]$, which directly quantifies the magnitude of step-wise changes, the covariance is less sensitive to small fluctuations.
If $\mathbb{E}[S(t)^2]$ is not already available from prior computations, evaluating $\mathbb{E}[(S(t+1) - S(t))^2]$ may require additional effort to compute it. However, in our model, $\mathbb{E}[S(t)^2]$ has already been used in previous moment equations when estimating the parameters $\pi_+$ and $\pi_-$. Therefore, using $\mathbb{E}[(S(t+1) - S(t))^2]$ not only avoids extra computational cost but also provides a more informative and sensitive measure of dynamic changes.

We chose the one--lag difference rather than a \( k \)--lag difference with \( k > 1 \), because as \( k \) increases, the dependence between two observations \( k \) steps apart decreases. One can consider the extreme case where \( k \) is very large; in that case, the observations at times \( t \) and \( t + k \) are nearly independent. In other words, the information contained in the quantity \( \mathbb{E}\big[(S(t+k) - S(t))^2\big] \) diminishes as \( k \) grows. Additionally, the computation of \( \mathbb{E}\big[(S(t+k) - S(t))^2\big] \) becomes increasingly complex with larger \( k \), making the trade-off unfavorable: we pay more computationally but gain less information. We also conducted numerical experiments for \( k = 2 \) (not presented in this paper), and the results yielded less accurate estimates, with larger variances and poorer agreement with the true values compared to the case \( k = 1 \).

We then define the empirical moment
\[
M_{3,K} := \frac{1}{K-1} \sum_{t = 1}^{K-1} \big(S(t+1) - S(t)\big)^2,
\]
and use the moment equation
\begin{equation} \label{eq:estimation2}
    \mathbb{E}\big[(S(t+1) - S(t))^2\big] = M_{3,K},
\end{equation}
to estimate the parameter $\alpha$. The overall estimation procedure consists of the following steps:
\begin{itemize}
    \item \textbf{Step 1}: Estimate $\pi_+$ and $\pi_-$, denoted by $\hat{\pi}_+$ and $\hat{\pi}_-$, by solving the moment equations in~\eqref{eq:estimation1};
    \item \textbf{Step 2}: Estimate $\alpha$, denoted by $\hat{\alpha}$, by solving~\eqref{eq:estimation2}, using the estimates $\hat{\pi}_+$ and $\hat{\pi}_-$ obtained in Step 1 in place of $\pi_+$ and $\pi_-$ in the expression of $\mathbb{E}\big[(S(t+1) - S(t))^2\big]$.
\end{itemize}

\section{Numerical examples}
\label{sec:simulations}

\begin{figure}[!h]
\centering
\subcaptionbox{$0.8986 \, (0.0078)$}
{\includegraphics[width=0.9\linewidth]{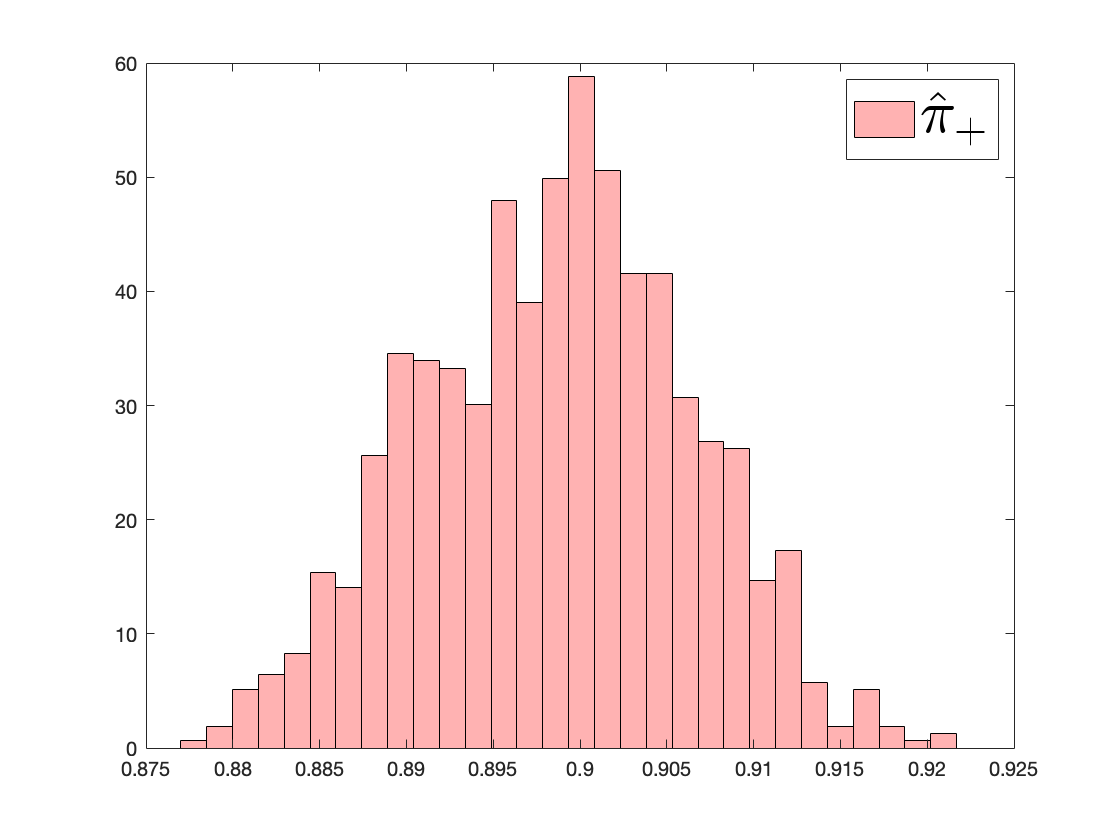}}
\subcaptionbox{$0.4015  \, (0.0079)$}
{\includegraphics[width=0.9\linewidth]{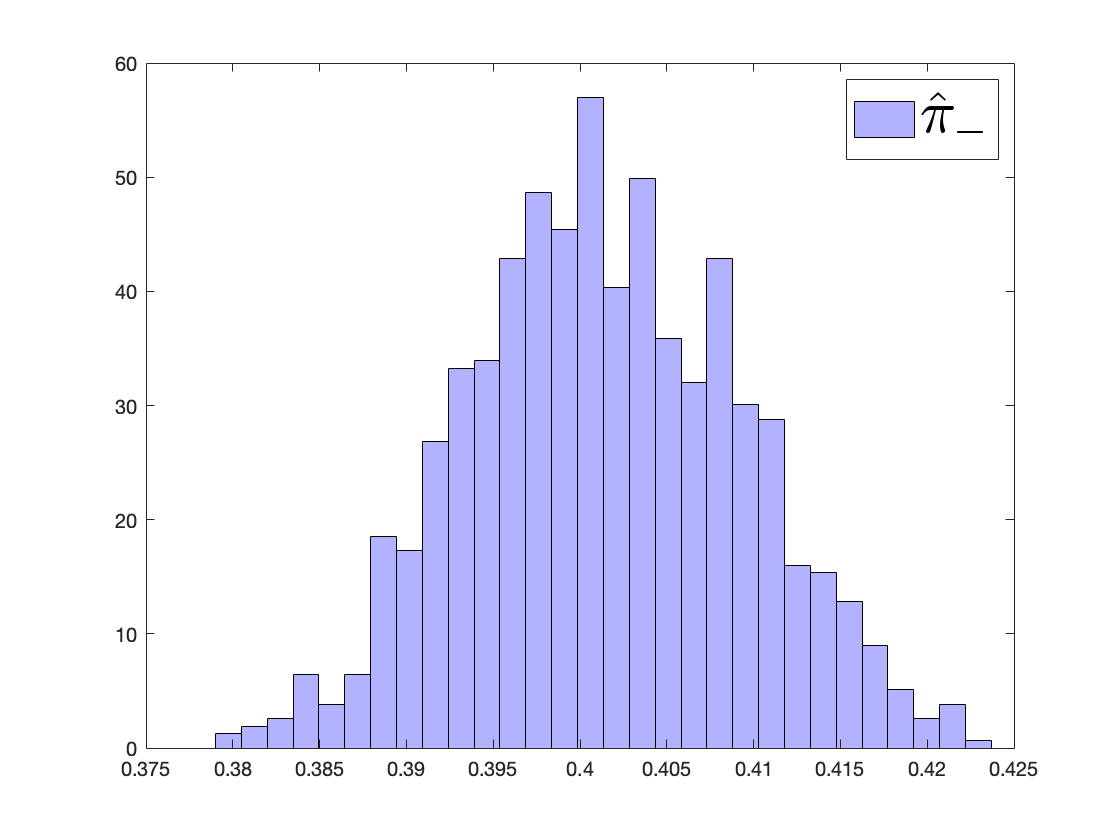}}
\subcaptionbox{$0.3003 \, (0.0072)$}
{\includegraphics[width=0.9\linewidth]{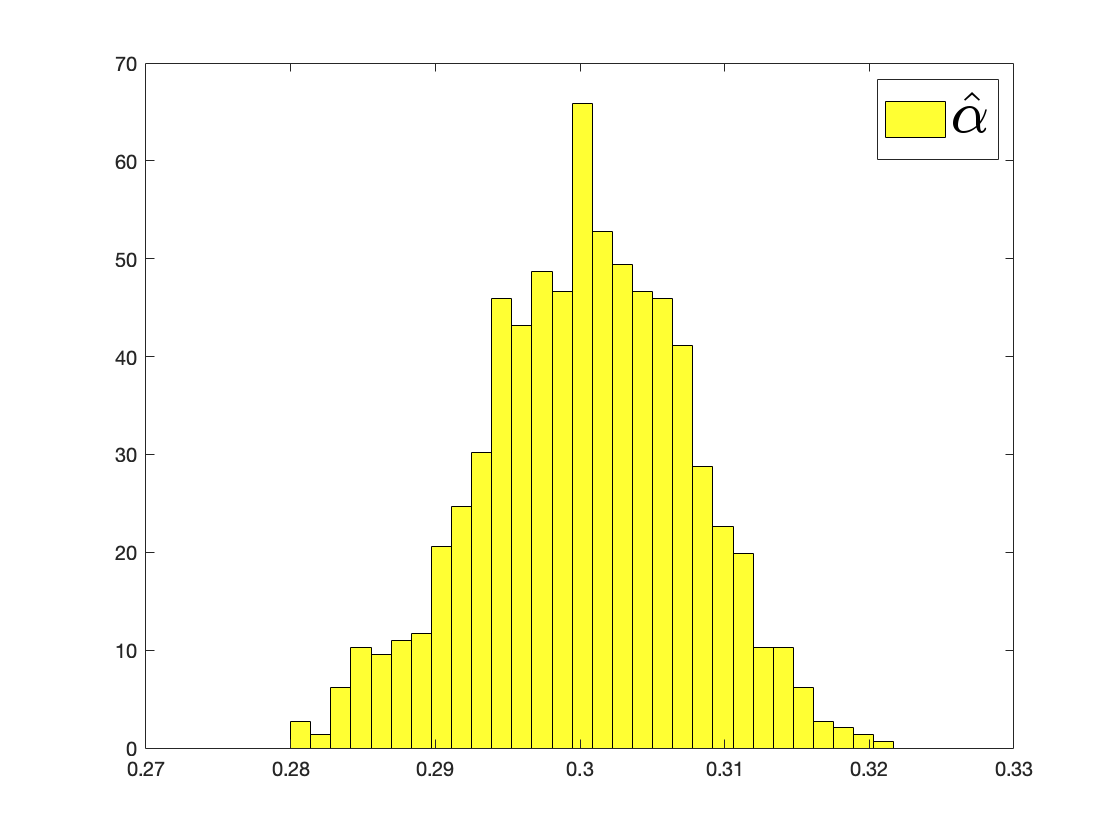}}
\caption{Estimation of the parameters $\pi_+$, $\pi_-$, $\alpha$ when $n = 10$, using the methodology described in Section~\ref{sec:methodology}, with $f(\pi_+, \pi_-) = (\pi_+ + \pi_-)/2$.} 
\label{fig:Estimation1}
\end{figure}

\begin{figure}[!h]
\centering
\subcaptionbox{$0.8996 \, (0.0080)$}
{\includegraphics[width=0.9\linewidth]{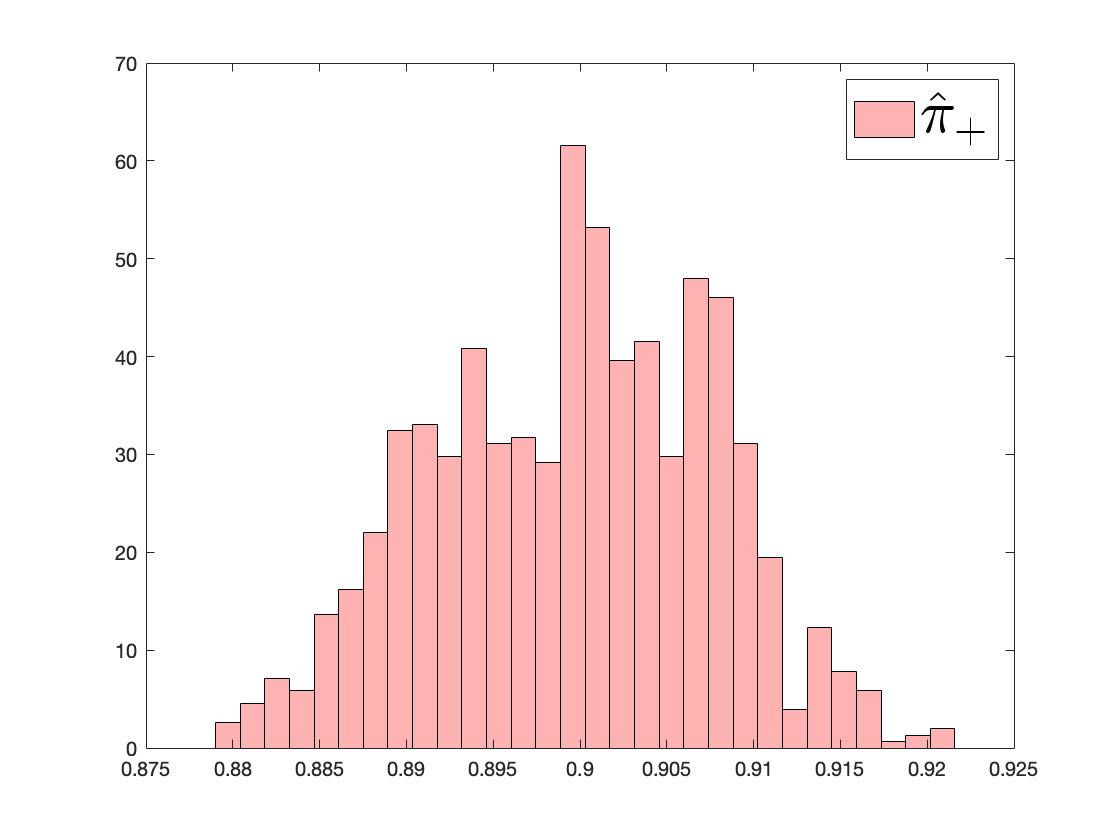}}
\subcaptionbox{$0.4003  \, (0.0040)$}
{\includegraphics[width=0.9\linewidth]{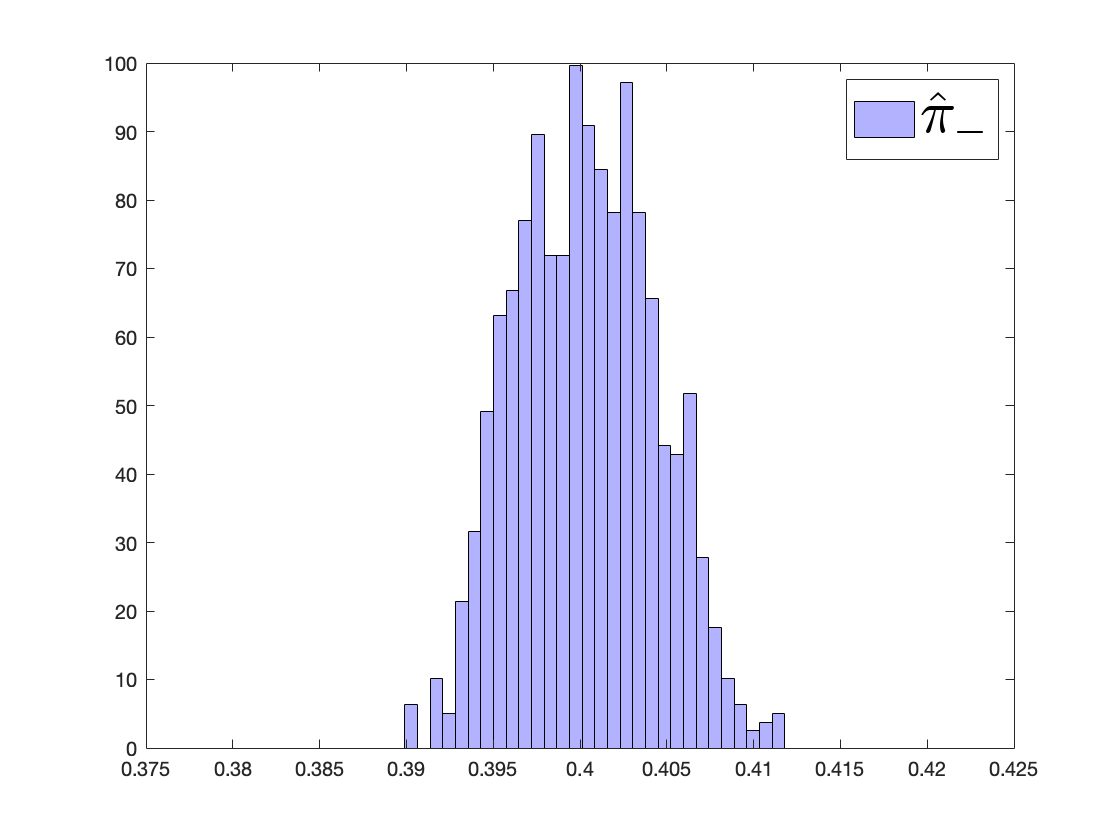}}
\subcaptionbox{$0.3001 \, (0.0093)$}
{\includegraphics[width=0.9\linewidth]{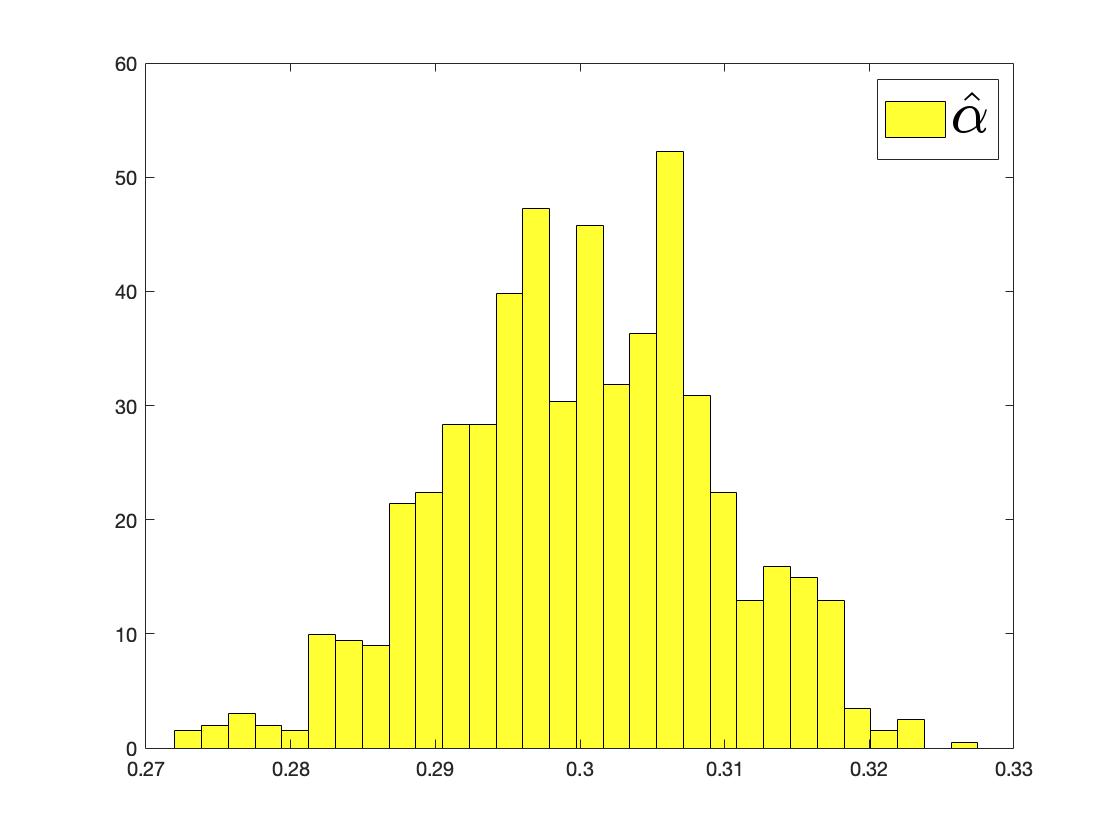}}
\caption{Estimation of the parameters $\pi_+$, $\pi_-$, $\alpha$ when $n = 10$, using the methodology described in Section~\ref{sec:methodology}, with $f(\pi_+,\pi_-)=\pi_+\pi_-/(\pi_++\pi_-)$.} 
\label{fig:Estimation2}
\end{figure}

We conduct numerical experiments with $L = 1000$ independent runs. In each run, we simulate the system for $K = 10^5$ observation points and estimate the parameters $\pi_+$, $\pi_-$ and $\alpha$. Let $\hat{\pi}_{+}^{(\ell)}$, $\hat{\pi}_{-}^{(\ell)}$, and $\hat{\alpha}^{(\ell)}$ 
denote the estimates obtained in the $\ell$-th run, for $\ell \in \{1, \ldots, L\}$. 
We define the sample mean and sample variance of the estimates of $\alpha$ as
\[
\bar{\alpha}^{(L)} = \displaystyle \dfrac{1}{L} \sum_{\ell=1}^L \hat{\alpha}^{(\ell)}
\]
and
\[
\sigma^2[\bar{\alpha}^{(L)}] = \displaystyle \dfrac{1}{L - 1} \sum_{\ell=1}^L (\hat{\alpha}^{(\ell)} - \bar{\alpha}^{(L)})^2,
\] 
respectively. Analogous definitions are used for the estimates of the other parameters $\pi_+$ and $\pi_-$. The experiments are conducted for two types of link functions $f(\pi_+, \pi_-)$, and the results are presented in Figures~\ref{fig:Estimation1} and~\ref{fig:Estimation2}. For each case, the sample mean and variance of the estimates are reported in the format $\bar{\alpha}_L \; (\sigma^2[\bar{\alpha}_L])$ displayed below the corresponding plots. Estimates for $\pi_+$ and $\pi_-$ are reported using the same notation. The true parameter values used in the simulations are:
\[
\pi_+ = 0.9, \quad \pi_- = 0.4, \quad \alpha = 0.3.
\]

When the link function is symmetric, the observations remain unchanged if \(\pi_+\) and \(\pi_-\) are swapped. Consequently, when estimating \(\pi_+\) and \(\pi_-\), it may occur that \(\widehat{\pi}_+\) is close to \(\pi_-\) and \(\widehat{\pi}_-\) is close to \(\pi_+\). To visualize the estimates in such cases, when \(\pi_+\) and \(\pi_-\) are known to be different and the link function is symmetric, we assign the larger estimate to \(\widehat{\pi}_+\) and the smaller estimate to \(\widehat{\pi}_-\). Note that this issue does not arise in the asymmetric case. We explain this in more details in \textbf{Asymmetric cases} below.

From Figures~\ref{fig:Estimation1} and~\ref{fig:Estimation2}, we observe that the estimates of the three parameters are quite accurate: the sample means are close to the true values, and the variances are within an acceptable range. The variances of the estimates for \( \pi_+ \) and \( \alpha \) are similar across both figures. However, the variance of the estimate for \( \pi_- \) is noticeably larger when the link function is given by \( f(\pi_+, \pi_-) = (\pi_+ + \pi_-)/2 \), compared to when it is \( f(\pi_+, \pi_-) = \pi_+ \pi_- / (\pi_+ + \pi_-) \). 
This difference may be attributed to the derivative of the link function with respect to \( \pi_- \). When the derivative is larger than \( 1/2 \), changes in \( \pi_- \) have a greater effect on the observed edge statistics, making \( \pi_- \) easier to infer from the data.

In what follows, we perform a sensitivity analysis to examine the robustness and potential bias of the estimator with respect to changes in the underlying parameters. Specifically, we assess the estimation performance for different values of $n$, $\pi_+$, $\pi_-$ and $\alpha$, including the case of asymmetric link function $f$ and of observations that are contaminated by noise.

\begin{figure}[!h]
\centering
\subcaptionbox{$0.8998 \, (0.0069)$}
{\includegraphics[width=0.9\linewidth]{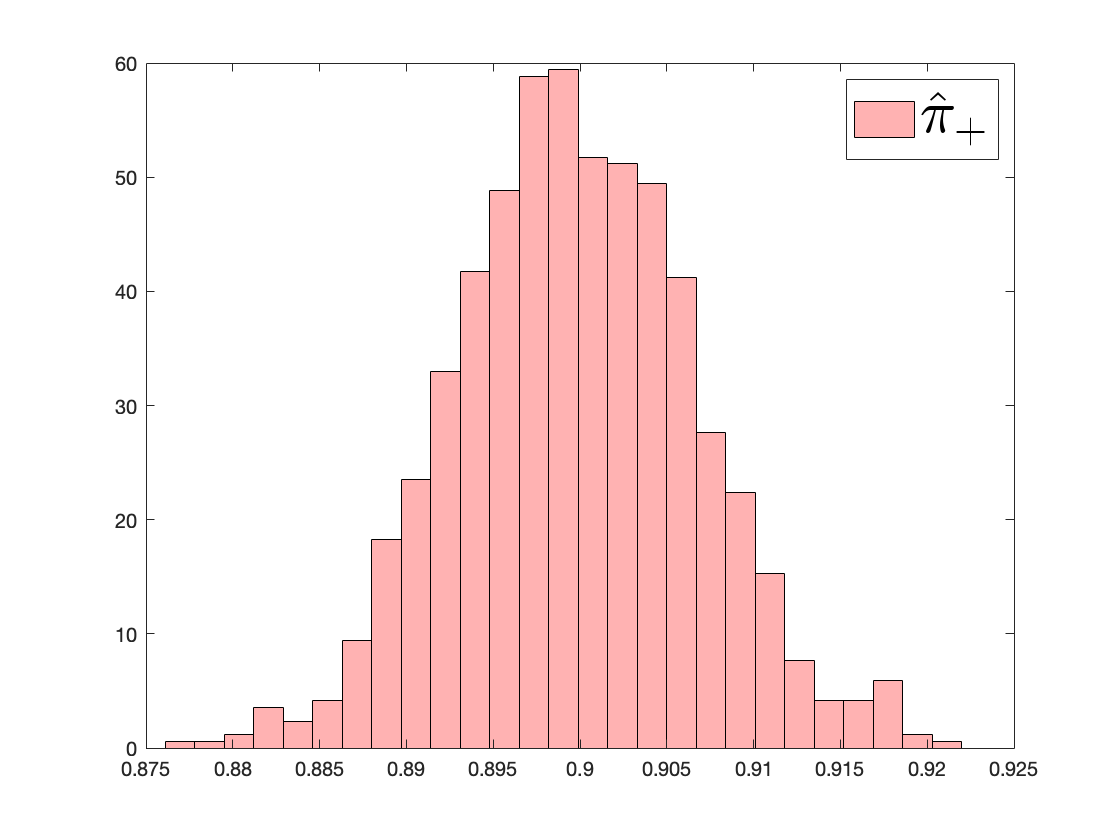}}
\subcaptionbox{$0.4003  \, (0.0027)$}
{\includegraphics[width=0.9\linewidth]{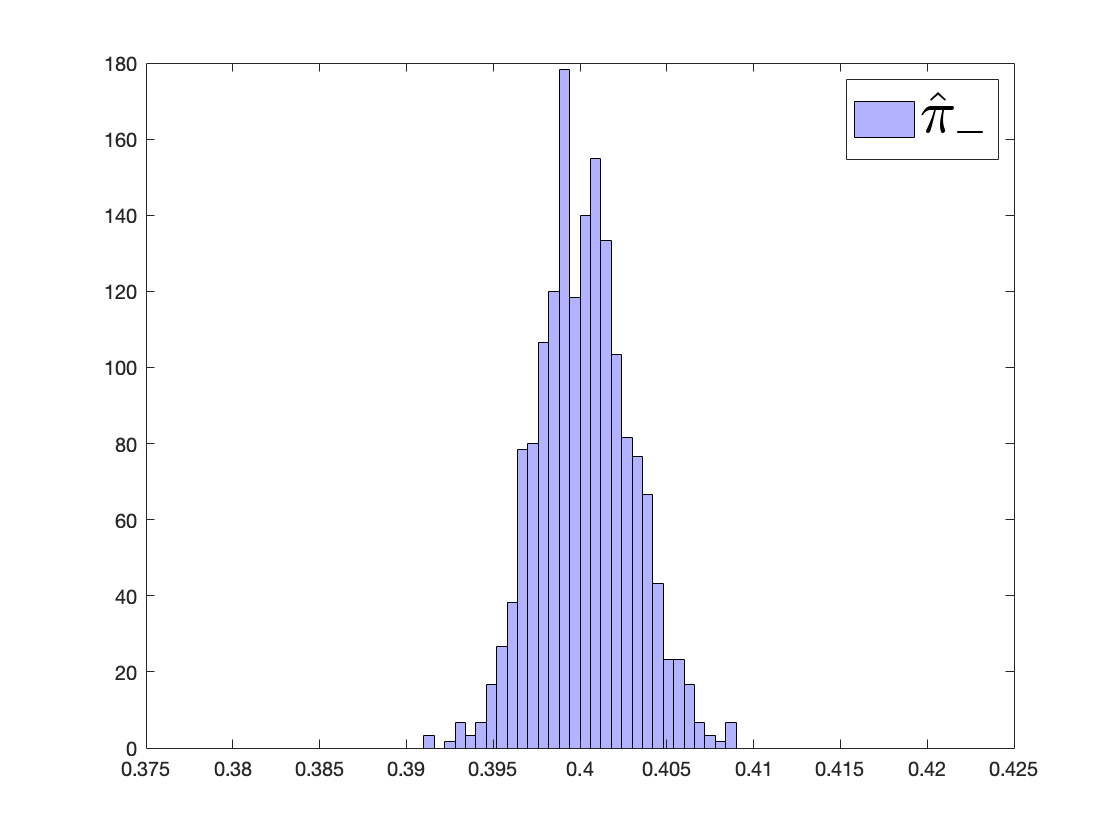}}
\subcaptionbox{$0.2997 \, (0.0087)$}
{\includegraphics[width=0.9\linewidth]{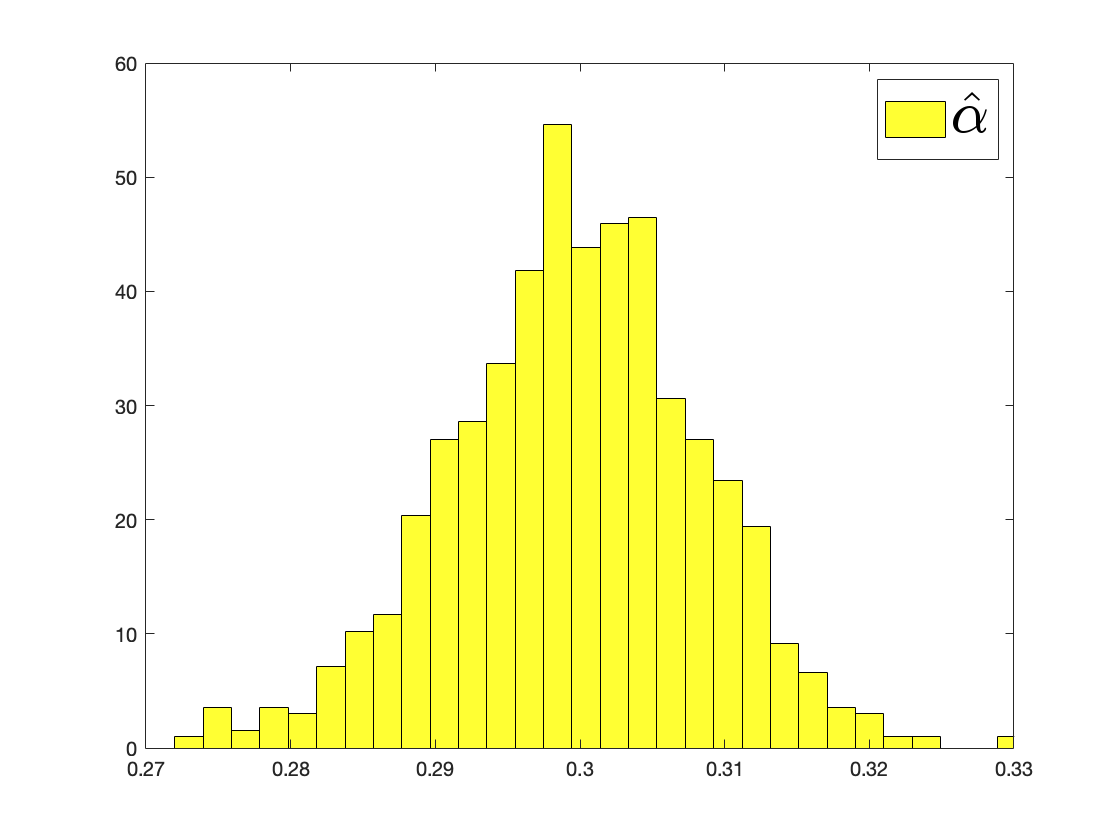}}
\caption{Estimation of the parameters $\pi_+$, $\pi_-$, $\alpha$ when $n = 5$, using the methodology described in Section~\ref{sec:methodology}, with $f(\pi_+,\pi_-)=\pi_+ \pi_- / (\pi_+ + \pi_-)$.} 
\label{fig:Estimation3}
\end{figure}

\begin{figure}[!h]
\centering
\subcaptionbox{$0.8985 \, (0.0084)$}
{\includegraphics[width=0.9\linewidth]{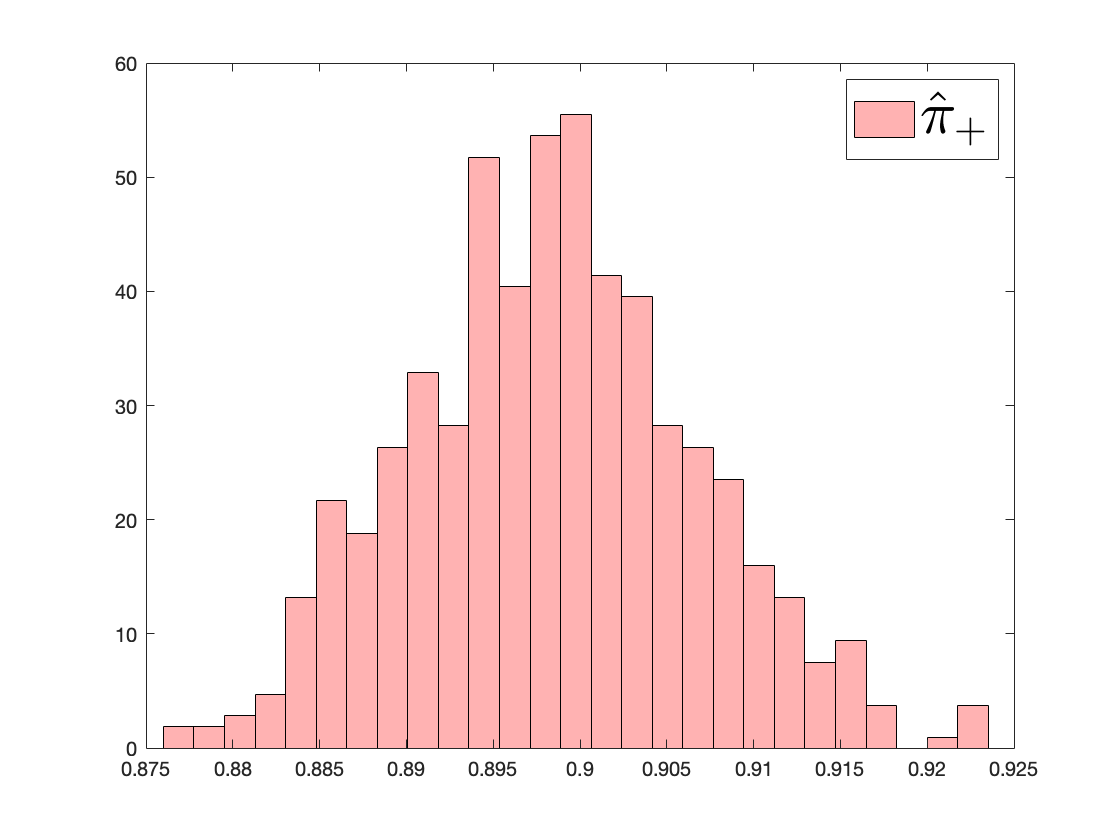}}
\subcaptionbox{$0.4007  \, (0.0041)$}
{\includegraphics[width=0.9\linewidth]{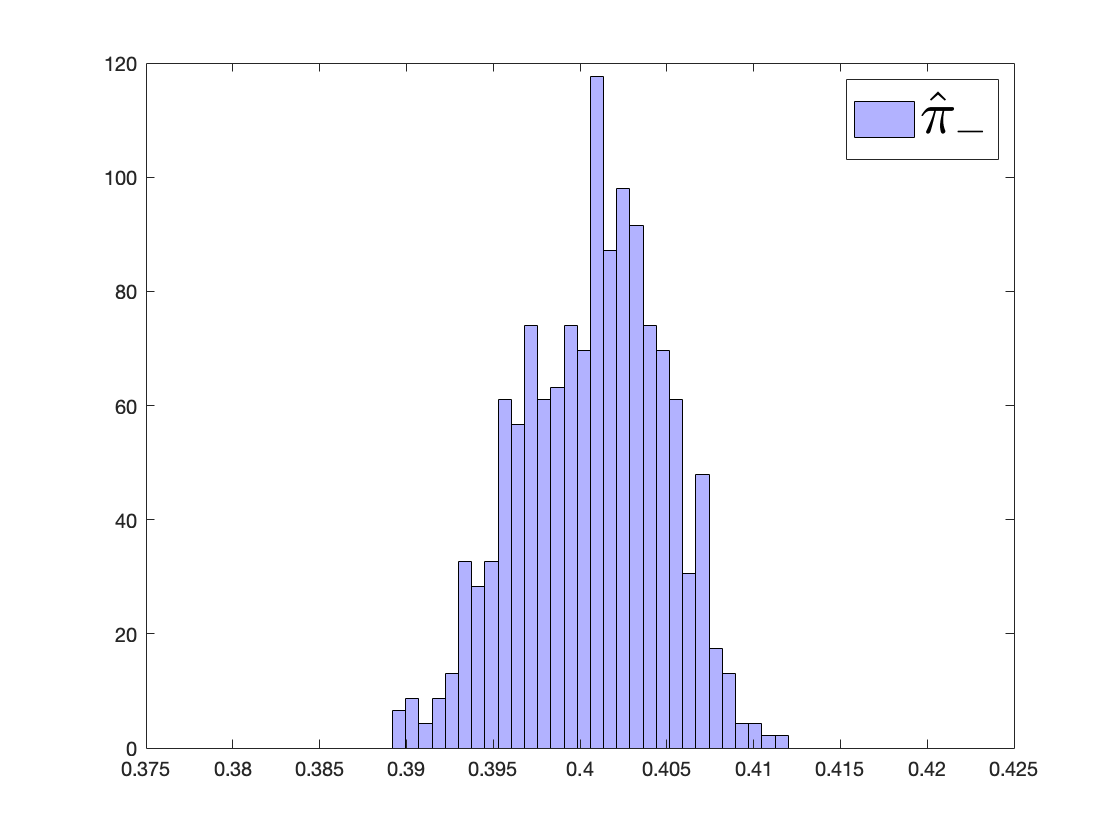}}
\subcaptionbox{$0.3010 \, (0.0098)$}
{\includegraphics[width=0.9\linewidth]{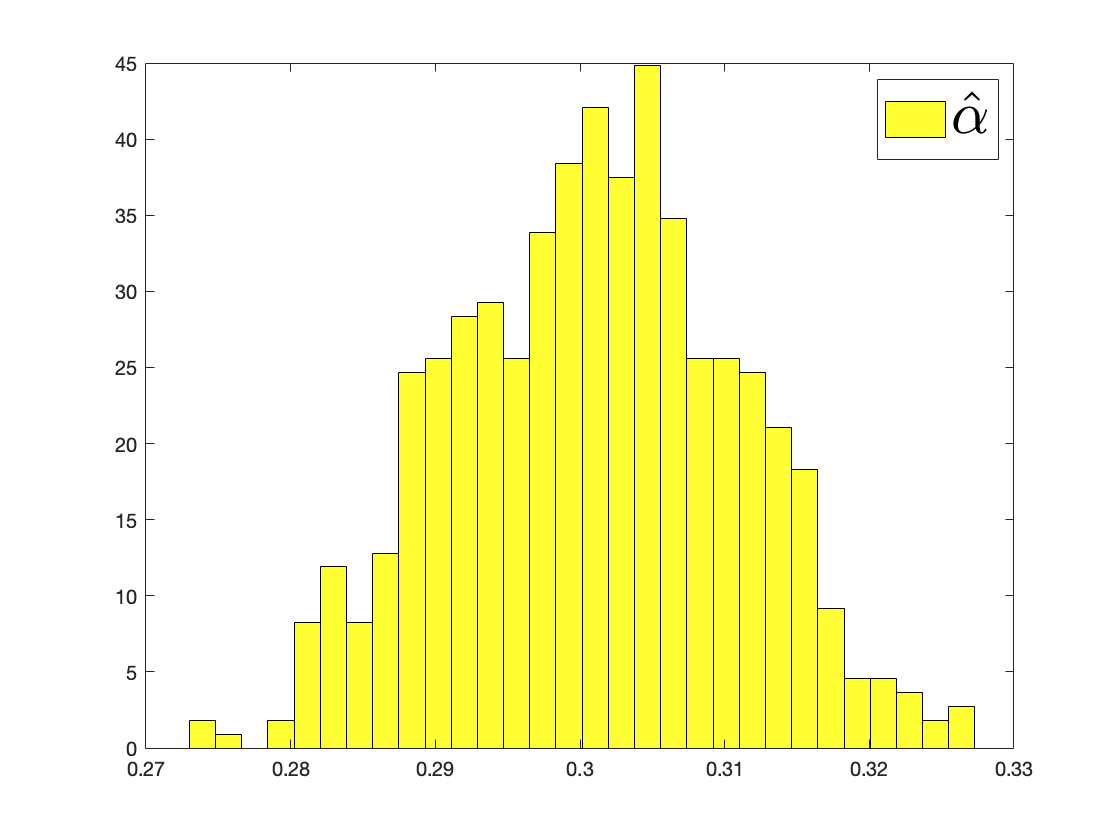}}
\caption{Estimation of the parameters $\pi_+$, $\pi_-$, $\alpha$ when $n = 20$, using the methodology described in Section~\ref{sec:methodology}, with $f(\pi_+,\pi_-)=\pi_+ \pi_- / (\pi_+ + \pi_-)$.} 
\label{fig:Estimation4}
\end{figure}

\medskip
\textbf{Different $n$.}
To examine the effect of the number of vertices \(n\) on the estimation, we fixed the parameters $\pi_+$, $\pi_-$ and $\alpha$, together with the link function $f$, and let $n$ assume different values. In particular, we set
\[
\pi_+ = 0.9, \quad \pi_- = 0.4, \quad \alpha = 0.3,
\]
with
\[
f(\pi_+, \pi_-) = \dfrac{\pi_+ \pi_-}{\pi_+ + \pi_-},
\]
and $n\in\{5,20\}$. The corresponding estimates are shown in Figures~\ref{fig:Estimation3} and~\ref{fig:Estimation4}. We can observe that as \(n\) increases, the variance grows slightly; however, this variation is negligible. The means remain close to the true values, and the distributions exhibit approximately bell--shaped patterns. 

\medskip
\textbf{Different $\alpha$.}
We next investigate the effect of \(\alpha \in \{0.1, 0.3, 0.5, 0.7\}\) on the estimates, while keeping
\[
\pi_+ = 0.9, \quad \pi_- = 0.4
\]
fixed and focusing on \(n = 20\).
Table~\ref{tab:alpha} presents the estimates of \(\pi_+\), \(\pi_-\), and \(\alpha\), reporting the mean values with variances shown in parentheses.
\begin{table}
\centering
\begin{tabular}{c c c c}
\toprule
\multirow{2}{*}{$\alpha$} & $\bar{\pi}_+^{(100)} $ & $\bar{\pi}_-^{(100)}$  & $\bar{\alpha}^{(100)}$ \\
& $(\sigma^2[\bar{\pi}_+^{(100)}])$ & $(\sigma^2[\bar{\pi}_-^{(100)}])$ & $(\sigma^2[\bar{\alpha}^{(100)}])$\\
\midrule\midrule
\multirow{2}{*}{0.1} & 0.8997 & 0.4000 & 0.0976 \\
& (0.0109) & (0.0057) & (0.0206) \\
\midrule
\multirow{2}{*}{0.3} & 0.8995 & 0.4006 & 0.3006 \\
& (0.0077) & (0.0043) & (0.0102) \\
\midrule
\multirow{2}{*}{0.5} & 0.8986 & 0.4007 & 0.4998 \\
& (0.0071) & (0.0036) & (0.0075) \\
\midrule
\multirow{2}{*}{0.7} & 0.8993 & 0.4005 & 0.6996 \\
& (0.0063) & (0.0036) & (0.0066)
\\
\bottomrule
\end{tabular}
\caption{For $n=20$ and $L =100$ runs, estimation of the parameters $\pi_+$, $\pi_-$, $\alpha$ for $\pi_+=0.9$, $\pi_-=0.4$ and $\alpha\in\{0.1,0.3,0.5,0.7\}$.}
\label{tab:alpha}
\end{table}

We can see that the means are all close to the true values, while the variances decrease as \(\alpha\) increases. When \(\alpha\) is small, the edge dynamics plays a more prominent role, leading to higher values of \(\mathrm{Var}[M_{1,K}]\) and \(\mathrm{Var}[M_{2,K}]\), even though \(\mathbb{E}[S(t)]\) and \(\mathbb{E}[S^2(t)]\) do not depend on \(\alpha\). Consequently, the increased variability results in higher variance of the estimates.

\medskip
\textbf{Different \(\mathbf \pi_+\) and \(\mathbf \pi_-\)}.
Next, we analyze the effect of \(\pi_+\) and \(\pi_-\) by considering cases where they are very different and where they are equal. Specifically, we set $n = 20$ and examine the cases
\[
\pi_+ = 0.9, \quad \pi_- = 0.2 \quad \text{and} \quad \pi_+ = \pi_- = 0.5.
\]
The results are presented in Table~\ref{tab:pis}. For the case \(\pi_+ = \pi_- = 0.5\), since the parameters are identical, we did not assign the smaller estimate to \(\pi_-\) and the larger estimate to \(\pi_+\) (doing so would introduce bias in the estimates).
\begin{table}
\centering
\begin{tabular}{c c c c c}
\toprule
\multirow{2}{*}{$\pi_+$}  & \multirow{2}{*}{$\pi_-$} & $\bar{\pi}_+^{(100)} $ & $\bar{\pi}_-^{(100)}$  & $\bar{\alpha}^{(100)}$ \\
 & & $(\sigma^2[\bar{\pi}_+^{(100)}])$ & $(\sigma^2[\bar{\pi}_-^{(100)}])$ & $(\sigma^2[\bar{\alpha}^{(100)}])$\\
\midrule\midrule
\multirow{2}{*}{0.9} & \multirow{2}{*}{0.2} & 0.9004 & 0.1991 & 0.2983 \\
 & & (0.0192) & (0.0203) & (0.0172) \\
\midrule
\multirow{2}{*}{0.5} & \multirow{2}{*}{0.5} & 0.5032 & 0.4967 & 0.2989 \\
 & & (0.0334) & (0.0334) & (0.0036) 
\\
\bottomrule
\end{tabular}
\caption{For \(n = 20\) and \(L = 100\) runs, estimation of the parameters \(\pi_+\), \(\pi_-\), \(\alpha\) for $(\pi_+,\pi_-)\in\{(0.9,0.2),(0.5,0.5)\}$ and $\alpha=0.3$.}
\label{tab:pis}
\end{table}

\medskip

\textbf{Noisy observations}.
We consider now the case in which the number of edges observed at each time point is not exactly \(S(t)\). Instead, the observed process is contaminated by a Gaussian white noise \(\varepsilon_t \sim \mathcal{N}(0,\sigma^2)\), which is independent of \(S(t)\). Thus, each observation can be written as $S(t) + \varepsilon_t$. Since the noise is independent and white,
\[
\mathbb{E}[S(t) + \varepsilon_t] = \mathbb{E}[S(t)], \ \ \mathbb{E}[(S(t) + \varepsilon_t)^2] = \mathbb{E}[S(t)^2] + \sigma^2,
\]
with covariance
\[
\mathrm{Cov}(S(t) + \varepsilon_t, \, S(r) + \varepsilon_r) = 
\begin{cases}
\mathrm{Cov}(S(t), S(r)), & t \neq r, \\
\mathrm{Var}(S(t)) + \sigma^2, & t = r.
\end{cases}
\]
which follows from
\[
\mathrm{Cov}(\varepsilon_t, \varepsilon_r) = 
\begin{cases}
\sigma^2, & t = r, \\
0, & t \neq r.
\end{cases}
\]
In other words, when $t \neq r$,
\[
\mathbb{E}[(S(t) + \varepsilon_t) (S(r) + \varepsilon_r)] = \mathbb{E}[S(t) \, S(r)],
\]  

To analyze the effect of noise, we first consider \(n = 5\) with \(\sigma = 0.1\): the results are shown in the first line of Table~\ref{tab:noise}. 
Increasing \(\sigma\) to 1 makes the estimation of \(\pi_+\) and \(\pi_-\) impossible, i.e., no parameter values can satisfy the moment equations. 
However, if we increase \(n\) to 20 while keeping \(\sigma = 1\), the parameters can still be estimated. 
This is because \(\mathbb{E}[S(t)^2]\) is small when \(n = 5\) and thus easily overwhelmed by the noise, whereas when \(n = 20\), \(\mathbb{E}[S(t)^2]\) is large, so the increase due to \(\sigma^2 = 1\) is negligible. 
Consequently, \(\pi_+\) and \(\pi_-\) remain well estimated.

Table~\ref{tab:noise} also presents the estimates 
for \((n,\sigma)=(20,1)\). 
For comparison, we also include the noise--free case, using the same seed 1989 for reproducibility. 
It can be seen that when \(n = 5\), the estimates are still fairly accurate, but the variances are somewhat larger under noise. 
When \(n = 20\), the noise has little effect on the estimates of \(\pi_+\) and \(\pi_-\), since the value of \(M_{2,k}\) is of order \(10^4\), making \(\sigma^2\) negligible in comparison. 

The estimator of $\alpha$ tends to be biased downward. Indeed, observe that
\[ 
\begin{array}{ll} 
\ee[((S(t+1)+\varepsilon_{t+1})-(S(t)+\varepsilon_{t}))^2] \\ 
\qquad \quad = 2 \big(\mathbb{E}[S(t)^2] + \sigma^2 - \mathbb{E}[S(t)S(t+1)]\big).
\end{array} 
\]
In our procedure, $\pi_+$ and $\pi_-$ are estimated first. These estimates correctly 
capture the contribution $\mathbb{E}[S(t)^2] + \sigma^2$, but they also inflate the value 
of $\mathbb{E}[S(t)S(t+1)]$. On the other hand, the mean of the empirical cross moments in the 
presence of noise remains unchanged. Since smaller values of $\alpha$ correspond to 
faster edge dynamics and therefore a smaller covariance $\mathbb{E}[S(t)S(t+1)]$, the 
model compensates by underestimating~$\alpha$.

\begin{table}
\centering
\begin{tabular}{c c c c c}
\toprule
\multirow{2}{*}{$n$} & \multirow{2}{*}{$\sigma$} & $\bar{\pi}_+^{(100)} $ & $\bar{\pi}_-^{(100)}$  & $\bar{\alpha}^{(100)}$ \\
& & $(\sigma^2[\bar{\pi}_+^{(100)}])$ & $(\sigma^2[\bar{\pi}_-^{(100)}])$ & $(\sigma^2[\bar{\alpha}^{(100)}])$\\
\midrule\midrule
\multirow{2}{*}{5} & \multirow{2}{*}{0.1} & 0.8997 & 0.4005 & 0.2929 \\
& & (0.0200) & (0.0219) & (0.0097) \\

\midrule
\multirow{2}{*}{5} & \multirow{2}{*}{0} & 0.8998 & 0.4003 & 0.2997 \\
    &            & (0.0069) & (0.0027) & (0.0087) \\
\midrule
\multirow{2}{*}{20} & \multirow{2}{*}{1} & 0.9015 & 0.3988 & 0.2614 \\
& & (0.0058) & (0.0061) & (0.0056) 
\\
\midrule
\multirow{2}{*}{20} & \multirow{2}{*}{0} & 0.8994 & 0.4000 & 0.2997 \\
& & (0.0081) & (0.0037) & (0.0095) 
\\
\bottomrule
\end{tabular}
\caption{For \(n \in \{5,20\}\) and \(L = 100\) runs, estimation of the parameters \(\pi_+\), \(\pi_-\), \(\alpha\) under noisy observations for $\pi_+=0.9$, $\pi_-=0.4$ and $\alpha=0.3$}.
\label{tab:noise}
\end{table}

\begin{table}
\centering
\begin{tabular}{c c c c}
\toprule
\multirow{2}{*}{$f(\pi_+, \pi_-)$} & $\bar{\pi}_+^{(100)} $ & $\bar{\pi}_-^{(100)}$  & $\bar{\alpha}^{(100)}$ \\
& $(\sigma^2[\bar{\pi}_+^{(100)}])$ & $(\sigma^2[\bar{\pi}_-^{(100)}])$ & $(\sigma^2[\bar{\alpha}^{(100)}])$ \\
\midrule\midrule
\multirow{2}{*}{$\frac{\pi_++\pi_-}{2}$} & 0.9005 & 0.3993 & 0.2988 \\
& (0.0094) & (0.0102) & (0.0072) \\
\midrule
\multirow{2}{*}{$\frac{\pi_+\pi_-}{\pi_++\pi_-}$} & 0.8998 & 0.4003 & 0.2997 \\
& (0.0069) & (0.0027) & (0.0087) \\
\midrule
\multirow{4}{*}{$\frac{\pi_++3\pi_-}{4}$} & 0.8988 & 0.4006 & 0.3004 \\
& (0.0051) & (0.0030) & (0.0053) \\
& 0.3373 & 0.7374 & 0.3000 \\
& (0.1930) & (0.1154) & (0.0091) \\
\midrule
\multirow{4}{*}{$\frac{\pi_+\pi_-}{\pi_++3\pi_-}$} & 0.8990 & 0.4003 & 0.3009 \\
& (0.0079) & (0.0034) & (0.0113)
\\
& 0.4445 & 0.9432 & 0.0433 \\
& (0.0033) & (0.0082) & (0.0493)
\\
\bottomrule
\end{tabular}
\caption{For \(n =5\) and \(L = 100\) runs, estimation of the parameters \(\pi_+\), \(\pi_-\), \(\alpha\) under different link functions for $\pi_+=0.9$, $\pi_-=0.4$ and $\alpha=0.3$.}
\label{tab:asymmetric}
\end{table}

\begin{figure}
\centering
{\includegraphics[width=0.9\linewidth]{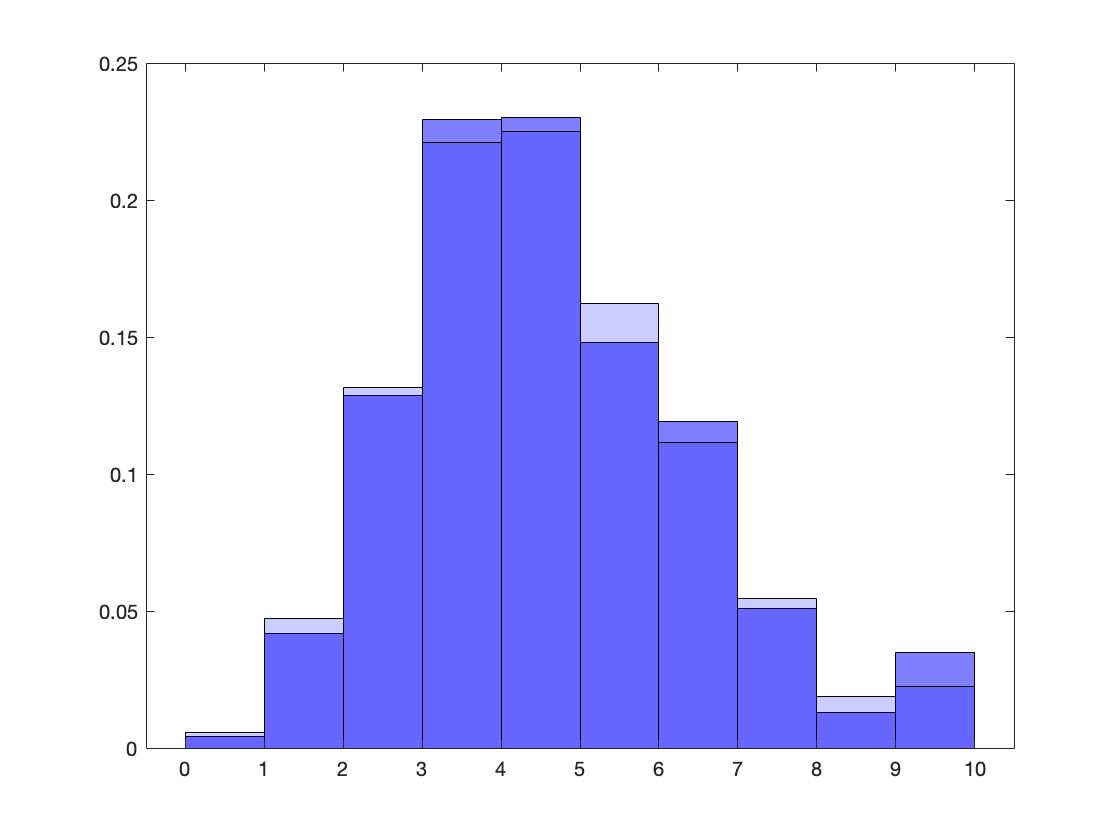}}
\caption{Histogram of the empirical distribution of $S(t)$ under the two different parameter sets $P_1$ (in dark blue) and $P_2$ (in light blue).}
\label{fig:EdgeStat}
\end{figure}

\medskip
\textbf{Asymmetric cases}.
As mentioned above, in the asymmetric case the estimates of \(\pi_+\) and \(\pi_-\) can fall into two possible sets, since both sets yield nearly the same values for \(M_{1,K}\) and \(M_{2,K}\). We present both sets in Table~\ref{tab:asymmetric}. 
To further analyze this situation, we fix $n = 5$, consider the link function $f(\pi_+, \pi_-) = \pi_+ \pi_-/(\pi_+ + 3 \pi_-)$, set the random seed to 1989, and run the simulation for \(K = 10^5\). We then compare the process \(\{S(t)\}\) under the following two parameter settings:
\[
P_1 : (\pi_+,\pi_-,\alpha)= (0.4445, \, 0.9432, \, 0.0433)
\]
and
\[
P_2 : (\pi_+,\pi_-,\alpha)=(0.8990, \, 0.4003, \, 0.3009).
\]
The histogram of the empirical distribution of \(\{S(t)\}\) is plotted in Figure \ref{fig:EdgeStat} under these two different choices of the parameters. We can observe that the two cases are {\it virtually} indistinguishable. This happens because the link function $f$ is asymmetric. Indeed, if $f$ is symmetric, the two parameters sets are just obtained by swapping the values of $\pi_+$ and $\pi_-$, showing that they correspond to the same edge dynamics.

\medskip
\textbf{Distinguish the links}.
With the same parameter set \((\pi_+, \pi_-, \alpha)\), different link functions are expected to yield different values of  
\[
\mathbb{E}[S(t)], \quad \mathbb{E}[S(t)^2], \quad \mathbb{E}\big[(S(t+1) - S(t))^2\big]. 
\]  
This implies that, by examining \(M_{1,K}\), \(M_{2,K}\), or \(M_{3,K}\), we should in principle be able to distinguish if the same link functions have been used or not. However, this is not always the case when considering only \(M_{3,K}\). We illustrate this below with a numerical example.

To make the comparison fair, we use Algorithm~\ref{alg:1} to generate \(M_{3,K}\) for different values of \(\alpha\), and present the results in Table~\ref{tab:distinguish}, where we compare the statistic \( M_{3,K} \), computed with \( K = 10^5 \) observations from \( L = 100 \) independent experiments, for two different choices of the link function \( f(\pi_+, \pi_-) \). In particular, we refer to this statistic as $M_{3,K}^1$ (resp.\ $M_{3,K}^2$) when $f(\pi_+,\pi_-)=(\pi_++\pi_-)/2$ (resp.\ $f(\pi_+,\pi_-)=(\pi_+\pi_-)/(\pi_++\pi_-)$).
As the parameter \( \alpha \) increases, updates predominantly occur at the vertices, while the edge configuration remains relatively stable. This results in strong temporal correlation in the edge states. In contrast, when \( \alpha \) is small, edges are resampled more frequently and independently at each step, leading to weaker but non-negligible correlation. This residual dependence arises because edge states are still sampled based on the current vertex configuration.
\begin{algorithm}[h!]
\SetAlgoLined
\caption{Compute $M_{3,K}$ using the \texttt{LinkFunction} for each $\alpha \in \text{par}$, where $\text{par} = \{0.1, 0.3, 0.5, 0.7\}$.}\label{alg:1}
Pick an $\alpha \in \text{par}$ and a link function;

Initialize empty list $\texttt{Record}^{\text{(LinkFunction)}}_\alpha$\;
$p_1 \gets 0.9$, $p_2 \gets 0.4$, $n \gets 20$\;
Set random seed to $1958$\;
\For{$\text{ind} \gets 1$ \KwTo $100$}{
    \textbf{Call} Simulation, and record $M_{3,K}$ in $\texttt{Record}^{\text{(LinkFunction)}}_\alpha$\; 
}
\end{algorithm}
\begin{table}
\centering
\begin{tabular}{cc c c c c}
\toprule
$n$ & $\alpha$ & $M_{3,K}^1 $ & $M_{3,K}^2$ & p-value & distinguishable \\
\hline
\hline
\multirow{2}{*}{10} & \multirow{2}{*}{0.1} & 16.3879 & 15.4057 & \multirow{2}{*}{$1.55\cdot10^{-45}$} & \multirow{2}{*}{Yes} \\
    & & (0.1070) & (0.1133) &  &  \\
    \midrule
    \multirow{2}{*}{10}
    & \multirow{2}{*}{0.3} & 13.7663 & 13.2127 & \multirow{2}{*}{$1.55\cdot10^{-45}$} & \multirow{2}{*}{Yes} \\
    & & (0.0930) & (0.0874) &  &  \\ 
    \midrule
    \multirow{2}{*}{10} & \multirow{2}{*}{0.5} & 10.9271 & 10.7625 & \multirow{2}{*}{$1.75\cdot10^{-19}$} & \multirow{2}{*}{Yes} \\
    & & (0.0902) & (0.0843) &  &  \\
    \midrule
    \multirow{2}{*}{10} & \multirow{2}{*}{0.6} & 9.3958 & 9.3599 & \multirow{2}{*}{0.0131} & \multirow{2}{*}{No} \\
    & & (0.0829) & (0.0836) &  &  \\ 
    \midrule
    \multirow{2}{*}{10} & \multirow{2}{*}{0.7} & 7.7074 & 7.7983 & \multirow{2}{*}{$1.06\cdot10^{-11}$} & \multirow{2}{*}{Yes} \\
    & & (0.0729) & (0.0807) &  &  \\
\hline
\hline
\multirow{2}{*}{20} 
    & \multirow{2}{*}{0.1} & 69.3639 & 64.9101 & \multirow{2}{*}{$1.55\cdot10^{-45}$} & \multirow{2}{*}{Yes} \\
    & & (0.4985) & (0.4209) &  &  \\ 
    \midrule
    \multirow{2}{*}{20} & \multirow{2}{*}{0.3} & 58.6061 & 55.7414 & \multirow{2}{*}{$1.55\cdot10^{-45}$} & \multirow{2}{*}{Yes} \\
    & & (0.4081) & (0.3889) &  &  \\ 
    \midrule
    \multirow{2}{*}{20} & \multirow{2}{*}{0.5} & 47.5214 & 46.0728 & \multirow{2}{*}{$1.23\cdot10^{-44}$} & \multirow{2}{*}{Yes} \\
    & & (0.3506) & (0.3425) &  &  \\ 
    \midrule
    \multirow{2}{*}{20} & \multirow{2}{*}{0.7} & 35.1759 & 34.9081 & \multirow{2}{*}{$2.45\cdot10^{-5}$} & \multirow{2}{*}{Yes} \\
    & & (0.3614) & (0.3859) &  &  \\ 
    \midrule
    \multirow{2}{*}{20} & \multirow{2}{*}{0.75} & 31.5881 & 31.5605 & \multirow{2}{*}{0.6766} & \multirow{2}{*}{No} \\
    & & (0.3737) & (0.3769) &  &  \\
\bottomrule
\end{tabular}
\caption{For $n \in\{ 10,20\}$, comparison of the statistic $M_{3,K}$, with $\pi_+ = 0.9$, $\pi_- = 0.4$, and $K=10^5$, corresponding to two different link functions $f$, for different values of $\alpha$.}
\label{tab:distinguish}
\end{table}

In the third and fourth columns of Table~\ref{tab:distinguish}, we applied MATLAB's \texttt{kstest2} function, a two-sample Kolmogorov--Smirnov test, whose null hypothesis is that the data come from the same continuous distribution. While the test generally rejects the null hypothesis at the 5\% significance level, it fails to do so at \(\alpha = 0.6\) for \(n = 10\) and at \(\alpha = 0.75\) for \(n = 20\), as shown in Table~\ref{tab:distinguish}.

\section{Conclusions}
\label{sec:conclusions}

The modeling and inference framework developed in this work is broadly applicable to real--world systems where the interaction between individual states and evolving network structure plays a central role. While the methodology is general, its relevance becomes particularly clear when applied to concrete domains.

In opinion dynamics, the framework can be used to study the evolution of social {\it polarization} or {\it consensus} formation in social networks. For instance, by fitting the model to empirical data from periods of low and high political activity, one could detect structural shifts in the interaction network, such as increasing segregation or echo chamber formation, based solely on aggregate link statistics. This makes it possible to track opinion fragmentation in real time, even when full knowledge of the network is unavailable.

In communication networks, the model provides a tool to infer protocol--level behavior from limited observations. For example, by training on edge counts in wireless sensor networks, one could monitor the responsiveness of the system to varying node activity. This enables diagnostics of performance degradation due to mobility, interference, or energy constraints. It could also be used to assess how resilient a protocol is under partial node failures or to detect emerging congestion patterns in decentralized architectures.

In brain networks, the model opens new possibilities for both clinical and experimental neuroscience. In the context of {\it disease detection}, training the model in healthy subjects allows estimation of baseline parameters that govern activity--dependent connectivity. Deviations in patient data, such as persistently low edge counts in motor or cognitive regions, can then be interpreted as evidence of neurons trapped in quiescent states, potentially supporting targeted therapeutic interventions. In {\it cognitive experiments}, the framework can be used to test hypotheses about attention or learning. For instance, whether focusing on a task increases the propensity for connectivity changes (i.e., modulates $\pi_+$ and $\pi_-$) in task--relevant areas.

The model was evaluated across a range of parameter settings, demonstrating the numerical consistency of the estimator. We further analyzed the results under varying modeling conditions. First, we introduced Gaussian noise into the observations to test the model’s robustness and examined its impact on the estimation performance. Second, to capture real--world heterogeneity, we allowed the influence of individual vertices to vary. For example, in opinion networks, some individuals may exert a stronger influence on others’ opinions, while others are less influential. To model this variability, we introduced both symmetric and asymmetric link functions and discussed potential issues of non--identifiability.

These examples illustrate how our inference method, despite relying on minimal observables, offers a flexible and data--efficient tool to extract meaningful insights across domains. Future extensions to fully coupled dynamics will further enhance its potential for applications in complex adaptive systems. We leave this as a natural and essential direction for future research.

\section*{Data availability}
The data that supports the findings of this article are openly available \cite{data}.

\appendix

\section{Proof of equation (\ref{eq:recurrence})}
\label{sec:appA}

    First, for any $0\leq k\leq n$, we will prove by induction over $k$ that
    \begin{equation}\label{eq:recurrence2}
    P_k=\binom{n}{k} P_0.
    \end{equation}
    The cases $k=0$ and $k=1$, since $P_1=nP_0$, are trivial to check. Suppose now that the claim holds for any $0\leq \ell\leq k$ and we want to prove it for $k+1$. By using \eqref{eq:recequation} together with the induction hypothesis, we get
    \[
    \begin{array}{ll}
   P_{k+1} &= \displaystyle\dfrac{n}{k+1} \binom{n}{k} P_0 - \dfrac{n-k+1}{k+1} \binom{n}{k-1} P_0  \\
   &= \displaystyle\dfrac{n!}{(k-1)! (n-k)!} \left( \dfrac{n}{k(k+1)} - \dfrac{1}{k+1} \right) P_0 \\
   &= \displaystyle\binom{n}{k+1} P_0.
   \end{array}
    \]
   This concludes the proof of \eqref{eq:recurrence2}. Finally, we prove \eqref{eq:recurrence} by showing that $P_0=1/2^n$ in order to have a probability density. This easily follows by the fact that $\sum_{k=0}^n \binom{n}{k}=2^n$.

\section{Proof of equation (\ref{eq:secondmoment})}
\label{sec:appB}
By \eqref{eq:defS}, note that we can write
\[
\begin{array}{ll}
S^2(t) &= \displaystyle\sum_{i<j} \mathbf{1}\{(i,j)\in E(t)\} \\
& \ + \displaystyle\sum_{i<j} \sum_{\substack{\ell<m \\ (i,j)\neq (\ell,m)}} \mathbf{1}\{(i,j)\in E(t), (\ell,m)\in E(t)\},
\end{array}
\]
which, at stationarity, leads to 
\begin{align*}
    &\ee[S^2(t)] - \ee[S(t)] \\
=   &\displaystyle\sum_{i<j} \sum_{\substack{\ell<m \\ (i,j)\neq (\ell,m)}} \pp\left( (i,j)\in E(t), (\ell,m)\in E(t) \right) \\
= &\displaystyle\sum_{\substack{i<j, \ \ell<m \\ (i,j)\neq (\ell,m)}} \sum_{k=0}^n\pp\left( (i,j),(\ell,m)\in E(t)| N_+(t)=k\right) P_k \\
= &\displaystyle \sum_{k=0}^n P_k \sum_{\substack{i<j, \ \ell<m \\ (i,j)\neq (\ell,m) \\ x_i(t)=+ \\ x_j(t)=+ }} \pp\left( (i,j),(\ell,m)\in E(t)| N_+(t)=k\right) 
\end{align*}
\begin{align*}
+&\displaystyle\sum_{k=0}^n P_k \sum_{\substack{i<j, \ \ell<m \\ (i,j)\neq (\ell,m) \\ x_i(t)=- \\ x_j(t)=-}} \pp\left( (i,j),(\ell,m)\in E(t)| N_+(t)=k\right) \\
+&\displaystyle\sum_{k=0}^n P_k \sum_{\substack{i<j, \ \ell<m \\ (i,j)\neq (\ell,m) \\ x_i(t)\neq x_j(t)}} \pp\left( (i,j),(\ell,m)\in E(t)| N_+(t)=k\right).
\end{align*}
We consider separately the three terms of the sum. The first one can be expanded as
\[
\begin{array}{ll}
\displaystyle \sum_{\substack{i<j, \ \ell<m \\ (i,j)\neq (\ell,m) \\ x_i(t)=+ \\ x_j(t)=+}} \pp\left( (i,j)\in E(t), (\ell,m)\in E(t)| N_+(t)=k\right) \\
  =\displaystyle \sum_{\substack{i<j \\ x_i(t)=x_j(t)=+}} \pp\left( (i,j)\in E(t)| N_+(t)=k\right)  \\
\quad \times \displaystyle\sum_{\substack{\ell<m \\ (i,j)\neq (\ell,m)}}\pp\left((\ell,m)\in E(t)| N_+(t)=k\right)  
 =\displaystyle \binom{k}{2}\pi_+ \\
\ \times \left\{ \left[ \binom{k}{2}-1 \right] \pi_+ + \binom{n-k}{2} \pi_- + k(n-k) f(\pi_+,\pi_-) \right\}.
\end{array}
\]
By arguing similarly for the other two terms, the claim follows.

\section{}
\label{sec:appC}
We show here that, if the link function $f$ is symmetric in the two variables and $(\pi_+,\pi_-)=(\pi_1,\pi_2)$ is a solution of the moment equations \eqref{eq:estimation1}, then also $(\pi_+,\pi_-)=(\pi_2,\pi_1)$ is. To this end, we will show that $\mathbb E[S_1(t)]= \mathbb E[S_2(t)]$ and $\mathbb E[S_1^2(t)]= \mathbb E[S_2^2(t)]$, where $S_1(t)$ (resp.\ $S_2(t)$) is the number of active edges at time $t$ for the model with $(\pi_+,\pi_-)=(\pi_1,\pi_2)$ (resp.\ $(\pi_+,\pi_-)=(\pi_2,\pi_1)$). First, for any $0\leq k \leq n$, by \eqref{eq:recurrence} we deduce that
\begin{equation}\label{eq:symbinomial}
    P_k = \dfrac{1}{2^n} \binom{n}{k} = \dfrac{1}{2^n} \binom{n}{n-k} = P_{n-k}.
\end{equation}
Thus, by using \eqref{eq:moment1} and \eqref{eq:symbinomial}, after using the change of variable $k'=n-k$, we deduce that $\mathbb E[S_1(t)]= \mathbb E[S_2(t)]$. After arguing in a similar way, by a direct computation we deduce that $\mathbb E[S_1^2(t)]= \mathbb E[S_2^2(t)]$ as well.

\end{document}